\documentclass[12pt]{article}
\usepackage{fullpage}
\usepackage{setspace}
\usepackage[colorlinks,linkcolor=blue, anchorcolor=red,citecolor=blue,CJKbookmarks=true]{hyperref}
\usepackage{epsfig, graphicx, floatrow, rotating, amssymb, bbm, amsthm, mathrsfs, dsfont, makeidx,color}
\usepackage{bm}
\usepackage{algorithmicx}
\usepackage{algpseudocode}
\usepackage[tbtags]{amsmath}
\usepackage[sectionbib]{natbib}
\newtheorem{theorem}{Theorem}[section]

\newtheorem{proposition}[theorem]{Proposition}

\newcommand{\ve}[1]{{\mbox{\boldmath ${#1}$}}}

\begin{document}

\title{Copula-based Partial Correlation Screening: a Joint and Robust Approach}
\author{ Xiaochao Xia$^{1, 2\ast}$ and Jialiang Li$^{2, 3, 4}$ \\
\textsl{ {\footnotesize $^{1}$College of Science, Huazhong Agricultural University, Wuhan, China} }  \\
\textsl{{ \footnotesize $^{2}$Department of Statistics and Applied Probability, National University of Singapore }} \\
\textsl{{ \footnotesize $^{3}$Duke-NUS Graduate Medical School }} \\
\textsl{{ \footnotesize $^{4}$ Singapore Eye Research Institute }} \\
}
\date{}
\maketitle

\begin{abstract}
Screening for ultrahigh dimensional features may encounter complicated issues such as outlying observations, heterogeneous or heavy-tailed distribution, multi-collinearity and confounding effects. Standard correlation-based marginal screening methods may be a weak solution to these issues. We contribute a novel robust joint screener to safeguard against outliers and distribution mis-specification for both the response variable and the covariates, and to account for external variables at the screening step. Specifically, we introduce a copula-based partial correlation (CPC) screener.  We show that the empirical process of the estimated CPC converges weakly to a Gaussian process and establish the sure screening property for CPC screener under very mild technical conditions, where we need not require any moment condition, weaker than existing alternatives in the literature. Moreover, our approach allows for a diverging number of conditional variables from the theoretical point of view. Extensive simulation studies and two data applications are included to illustrate our proposal.              
\end{abstract}

\section{Introduction}
\indent 
With the arrival of a big data era, ultrahigh dimensional data has become readily available from many business and scientific research fields, including medicine, genetics, finance and economics. Such massive data usually carry two common features: (i) the number of predictors or features can be tremendous and diverge to infinity with the sample size and (ii) the data distribution is very likely to be heteroscedastic and heavy-tailed for both the response and the covariates. These two features are observed in the two real data sets investigated in this paper. It is usually hoped that a variable screener can identify important predictors among numerous candidates. We note that eventually data scientists for such large scale data still need to construct a comprehensive model to accurately predict the future outcome. Thus a purely marginal screening approach as is usually adopted in the literature may not adequately serve the ultimate model building purpose.  We contribute a new screening method that addresses the above issues and complements the existing methodology.

Variable screening serves as a fast and efficient computing device. Abundant feature screening methods are proposed in recent decades, including the sure independence screening (SIS) by \cite{FL2008} who first established the sure screening property under Gaussian linear model, the sure independent ranking screening (SIRS, \cite{ZLLZ2011}), the Kendall's $\tau$ based screening (Kendall-SIS, \cite{LPZZ2012}), the distance correlation based screening (DC-SIS, \cite{LZZ2012}), the quantile-adaptive screening (QaSIS, \cite{HWH2013}), empirical likelihood screening (\cite{CTW2013,CTW2016}), the censored rank independence screening for lifetime data (CRIS, \cite{SLMJ2014}),
the screening method based on quantile correlation (QC-SIS, \cite{LLT2015}), the conditional quantile screening (CQ-SIS, \cite{WY2015}), the survival impaction index screening (SII, \cite{LZPH2016}),
the nonparametric independence screening (NIS, \cite{FFS2011}, \cite{CHLP2014,XYLZ2016}), among many others. These screening tools might suffer the following two drawbacks: First, almost all methods evaluate a marginal association between the response and the predictors without adjusting external variables. Therefore some jointly important markers may be incorrectly screened out if their marginal signal is not as strong as the spurious markers in the ranked list. On the other hand, marginally important variables may be jointly ineffective and hence including them in a multivariate model may lead to less convincing prediction (see eg. \cite{XJLZ2016}). To take into account the joint effects, a marginal feature screening is usually followed by an iterative calculation, such as the iterative SIS (ISIS) in \cite{FL2008}, which is computationally expensive and does not come with any theoretical guarantee. Secondly, distribution of the response and the covariates may be rather different from the light-tailed symmetric normal distribution and very often there are outliers affecting the computed screening indices. Some of the aforementioned procedures address the robustness of the response but to our knowledge none of the existing work addresses the robustness of the covariates yet, which is a harder problem with higher dimension.

We aim to tackle the two problems with a new screener. Specifically, to address the first issue, we develop a joint feature screening method by incorporating additional information.  Recently, a few conditional feature screening methods have been proposed. For instance, \cite{LLW2014} considered a sure independence screening procedure via conditional Pearson correlation coefficient through a kernel smoothing. Their method can be employed to handle ultrahigh dimensional varying-coefficient feature variables, which were investigated in \cite{FMD2014} and \cite{CHLP2014} as well. In addition, \cite{XLF2018} considered a robust screening method based on conditional quantile correlation, a generalised conception of \cite{LLT2015}. However, these authors only considered a single conditional variable. To extend to multivariate conditional variables, \cite{CLR2016} studied several confounding variables. \cite{BFV2016} extended \cite{FS2010})'s approach to allow for a portion of predictors as conditional variables. Our work provides a more general framework where all the ultrahigh dimensional markers and other low dimensional confounders can be jointly considered during the screening process. For the second issue, we incorporate robust copula-based correlation and partial correlation in our screening methods. The nonparametric copula is a well-known distribution-free summary measure and naturally leads to a screener robust against outliers and distribution mis-specification. To the best of our knowledge, there are very few works applying this classical dependence concept in high-dimensional setting.  \cite{XLF2018} proposed a robust conditional feature screening approach, however, their method performs only robustly against the response but not against the covariates. Another relevant recent work is \cite{MLT2017}.

The contribution of this paper can be summarised as follows. Firstly, we propose a doubly robust copula-based correlation (CC). Copula is a very popular bivariate function to model the nonlinear dependence between paired variates. See \cite{N2007} for an introduction to the copula. The CC characterises the empirical dependence between two random variables evaluated at a level pair and is invariant under monotone transformation for both variables. We study the asymptotic process properties of the CC. A marginal variable screening approach via CC (CC-SIS) can be performed and achieves the desired sure screening consistency (\cite{FL2008}). Secondly, extending copula-based correlation to copula-based partial correlation (CPC), we then construct a more general framework for joint screening. The importance of each marker is evaluated in the presence of conditional variables. This provides a fast way for conditional feature screening. CPC is also robust by its construction from a nonparametric estimation and thus may be more reliable than a similar approach in \cite{MLT2017} with a broader range of application. We provide both theoretical and numerical support for the proposed screening method. Our data analysis indicates that the final multivariate regression models built after our screening approach indeed predict the outcome with improved accuracy.

The rest of the paper is organised as follows. Section 2 presents CC, its empirical estimate and CC-SIS. Asymptotic properties for the estimated CC functions are established in Section 2 as well. Methodologies and large sample properties for the CPC and the CPC-SIS are presented in Section 3. Further implementation details on CPC-SIS for different cases are given in Section 4. Simulation studies and two applications are carried out in Section 5. Section 6 concludes the paper. All the technical proofs are relegated to the Appendix.

\section{Copula-based Correlation and Variable Screening}
 Consider two continuous random variables $X$ and $Y$. Let $F_{X}$ be the cumulative distribution function (CDF) of $X$, which is assumed to be right continuous. $F_{X}^{-1}(\tau)=\inf\{x: F_{X}(x)\geq \tau\}$ is the $\tau$ quantile of $F_{X}$. $F_{Y,X}$ is the joint CDF of $Y$ and $X$, and $F_{Y|X}$ is the conditional distribution function of $Y$ given $X$ with a density $f_{Y|X}$. We use $F_{n,X}$, $F_{n,X}^{-1}$ and $F_{n,Y,X}$ to denote empirical versions of $F_{X}$, $F_{X}^{-1}$ and $F_{Y,X}$, respectively, based on a sample of size $n$. Let $D[a, b]$ be the Banach space of all c\`{a}dl\`{a}g functions $z:[a,b] \mapsto \mathbb{R}$ on an interval $[a,b]\subset \bar{\mathbb{R}}$ equipped with the uniform norm, and $\ell^{\infty}([a,b]^2)$ denotes the collection of all bounded functions $z: [a,b]^2 \mapsto \mathbb{R}$. We use $\stackrel{d}{\to}$ to denote convergence in distribution.

\subsection{Copula-based Correlation}
We propose the following copula-based correlation (CC)
\begin{eqnarray}
  \varrho_{Y, X}(\tau, \iota) &=&  \frac{F_{Y,X}(F_{Y}^{-1}(\tau), F_{X}^{-1}(\iota)) -\tau\iota }{ \sqrt{\tau(1-\tau)\iota(1-\iota)}}, \qquad 0\le \tau,\iota \le 1,
    \label{eq:202}
\end{eqnarray}
where the first term in the numerator is a copula function $C(u,v)=F_{U,V}(u, v)$ with $U=F_{Y}(Y)$ and $V=F_{X}(X)$, evaluated at $(u,v)=(\tau, \iota)$ (see Corollary 2.3.7 of \cite{N2007}, p.22).  By a simple algebra, we have $ \varrho_{Y, X}(\tau, \iota) = \frac{\mathrm{E}[\psi_{\tau}(Y-F_{Y}^{-1}(\tau))\psi_{\iota}(X-F_{X}^{-1}(\iota))] }{\sqrt{\tau(1-\tau)\iota(1-\iota)}} = \frac{\mathrm{cov}(\psi_{\tau}(Y-F_{Y}^{-1}(\tau)), \psi_{\iota}(X-F_{X}^{-1}(\iota))) }{\sqrt{\tau(1-\tau)\iota(1-\iota)}}$,
where $\psi_{\tau}(u)=\tau-I(u \leq 0)$ and $I(\cdot)$ is the indicator function. Since  $\mathrm{var}(\psi_{\iota}(X-F_{X}^{-1}(\iota)))=\iota (1-\iota)$ and $\mathrm{var}(\psi_{\tau}(Y-F_{Y}^{-1}(\tau)))=\tau (1-\tau)$, $\varrho_{Y, X}(\tau, \iota)$ in (\ref{eq:202}) is indeed a legitimate correlation coefficient that lives between $-1$ and $1$. Like other correlation measures, CC equals 0 if $X$ and $Y$ are independent.

CC can measure the nonlinear dependence between $X$ and $Y$, thus incorporates all kinds of bivariate joint distribution of $X$ and $Y$. In addition, because the indicator function is unaffected by outliers and extreme values, CC is robust for certain heavy-tailed distribution for both $Y$ and $X$. We note that monotone transformation of $X$ and $Y$ does not alter the value of CC.

Given a sample of i.i.d. observations $\{(X_i, Y_i), i=1,\ldots, n\}$, we can construct an empirical estimate of $\varrho_{Y, X}(\tau, \iota)$ as
\begin{eqnarray}
  \widehat{\varrho}_{Y,X}(\tau, \iota)    =  \frac{F_{n, Y,X}(F_{n,Y}^{-1}(\tau), F_{n,X}^{-1}(\iota)) -\tau\iota }{ \sqrt{\tau(1-\tau)\iota(1-\iota)}}   
   = \frac{n^{-1}\sum_{i=1}^{n}\psi_{\tau}(Y_i-F_{n, Y}^{-1}(\tau))\psi_{\iota}(X_{i}-F_{n, X}^{-1}(\iota))  }{ \sqrt{\tau(1-\tau)\iota(1-\iota)}}
. \label{eq:202b}
\end{eqnarray}
Let $\sigma_{Y,X}(\tau, \iota)= F_{Y,X}(F_{Y}^{-1}(\tau), F_{X}^{-1}(\iota))$,
$\sigma_{X|Y}(\tau,\iota) =F_{X|Y=F_{Y}^{-1}(\tau)}(F_{X}^{-1}(\iota))$ and
$\sigma_{Y|X}(\tau,\iota) = F_{Y|X=F_{X}^{-1}(\iota)}(F_{Y}^{-1}(\tau))$.
In the following, we fix the level at $(\tau,\iota)$ and write $\sigma_{Y,X}$, $\sigma_{X|Y}$ and $\sigma_{Y|X}$ for simplicity.
Furthermore, define
\begin{eqnarray*}
  \xi(Y, X; \tau, \iota) &=& \frac{1}{\sqrt{\tau(1-\tau) \iota(1-\iota)}} \Big[ I(Y \leq F_{Y}^{-1}(\tau), X\leq F_{X}^{-1}(\iota)) \\
&& \quad -  \sigma_{X|Y}(\tau,\iota) I(Y\leq F_{Y}^{-1}(\tau))  -   \sigma_{Y|X}(\tau,\iota)  I(X\leq F_{X}^{-1}(\iota)) \Big].
\end{eqnarray*}
We have the weak convergence result for $\widehat{\varrho}_{Y, X}(\tau, \iota)$ established in the next theorem.
\begin{theorem}  \label{th1}
Let $0<a<b<1$  and suppose that marginal distributions $F_{X}$ and $F_{Y}$ are continuously differentiable on the intervals $[F_{X}^{-1}(a)-\varepsilon, F_{X}^{-1}(b)+\varepsilon]$ and $[F_{Y}^{-1}(a)-\varepsilon, F_{Y}^{-1}(b)+\varepsilon]$ with positive derivatives $f_{X}$ and $f_{Y}$, respectively, for some $\varepsilon >0$. Furthermore, assume that conditional density functions $f_{Y|X}$ and $f_{X|Y}$ are continuous on the product of these intervals. Then
\begin{eqnarray*}
   \sqrt{n}\{ \widehat{\varrho}_{Y, X }(\tau, \iota) -  \varrho_{Y, X}(\tau, \iota) \}  \stackrel{w}{\rightsquigarrow }    \mathbb{G}_{Y,X}(\tau, \iota)
\end{eqnarray*}
in $\ell^{\infty}([a,b]^2)$, where $\stackrel{w}{\rightsquigarrow}$ denotes "converge weakly", and $\mathbb{G}_{Y,X}(\tau, \iota)$  is Gaussian process with mean zero and covariance function
$\Omega_{1}(\tau_1,\iota_1;\tau_2,\iota_2) \equiv \mathrm{E}\{[\xi(Y,X; \tau_1, \iota_1) - \mathrm{E} \xi(Y, X; \tau_1, \iota_1)]\times [\xi(Y, X; \tau_2, \iota_2) - \mathrm{E} \xi(Y,X; \tau_2, \iota_2)]\}$.
\end{theorem}
We may explicitly write the covariance function $ \Omega_{1}(\tau_1,\iota_1;\tau_2,\iota_2)=$
\begin{eqnarray*}
&&  \big\{ F_{Y,X}(F_{Y}^{-1}(\tau_1\wedge \tau_2 ), F_{X}^{-1}(\iota_1 \wedge \iota_2))-F_{Y,X}(F_{Y}^{-1}(\tau_1), F_{X}^{-1}(\iota_1))F_{Y,X}(F_{Y}^{-1}(\tau_2), F_{X}^{-1}(\iota_2)) \\
&& -\sigma_{X|Y}(\tau_2,\iota_2)[ F_{Y,X}(F_{Y}^{-1}(\tau_1\wedge \tau_2 ), F_{X}^{-1}(\iota_1)) -  F_{Y,X}(F_{Y}^{-1}(\tau_1 ), F_{X}^{-1}(\iota_1))\tau_2 ]   \\
&& -\sigma_{Y|X}(\tau_2,\iota_2)[F_{Y,X}(F_{Y}^{-1}(\tau_1), F_{X}^{-1}(\iota_1 \wedge \iota_2)) - F_{Y,X}(F_{Y}^{-1}(\tau_1), F_{X}^{-1}(\iota_1))\iota_2]  \\
&& - \sigma_{X|Y}(\tau_1, \iota_1)[F_{Y,X}(F_{Y}^{-1}(\tau_1 ), F_{X}^{-1}(\iota_2)) - F_{Y,X}(F_{Y}^{-1}(\tau_2 ), F_{X}^{-1}(\iota_2))\tau_1] \\
&& + \sigma_{X|Y}(\tau_1,\iota_1)\sigma_{X|Y}(\tau_2,\iota_2)(\tau_1 \wedge \tau_2 -\tau_1\tau_2)  \\
&& + \sigma_{X|Y}(\tau_1,\iota_1)\sigma_{Y|X}(\tau_2,\iota_2)[F_{Y,X}(F_{Y}^{-1}(\tau_1 ), F_{X}^{-1}(\iota_2))-\tau_1\iota_2]  \\
&& - \sigma_{Y|X}(\tau_1,\iota_1)[ F_{Y,X}(F_{Y}^{-1}(\tau_2), F_{X}^{-1}(\iota_1 \wedge \iota_2)) -F_{Y,X}(F_{Y}^{-1}(\tau_2 ), F_{X}^{-1}(\iota_2))\iota_1  ] \\
&& + \sigma_{Y|X}(\tau_1,\iota_1) \sigma_{X|Y}(\tau_2,\iota_2) [F_{Y,X}(F_{Y}^{-1}(\tau_2), F_{X}^{-1}(\iota_1)) -\tau_2\iota_1 ]  \\
&& + \sigma_{Y|X}(\tau_1,\iota_1) \sigma_{Y|X}(\tau_2,\iota_2)(\iota_1 \wedge \iota_2 - \iota_1\iota_2) \big\}/[\tau_1(1-\tau_1)\iota_1(1-\iota_1)\tau_2(1-\tau_2)\iota_2(1-\iota_2)]^{1/2}.
\end{eqnarray*}
In particular, at fixed $(\tau, \iota)$, if $\varrho_{Y, X}(\tau, \iota)=0$, then $\sqrt{n} \widehat{\varrho}_{Y, X }(\tau, \iota) \stackrel{d}{\to }  N(0, \Omega_{1})$, where
$ \Omega_{1}\equiv \Omega_{1}(\tau,\iota; \tau,\iota)  =  \big\{ \sigma_{Y,X} - \sigma_{Y,X}^2   + (\tau -\tau^2)\sigma_{X|Y}^2 + ( \iota -\iota^2) \sigma_{Y|X}^2
 - 2(1-\tau)\sigma_{Y,X} \sigma_{X|Y} - 2(1- \iota) \sigma_{Y,X}\sigma_{Y|X}
  + 2[\sigma_{Y,X} - \tau\iota]\sigma_{X|Y} \sigma_{Y|X}
  \big\} /[\tau(1-\tau)\iota(1-\iota)]$.
If $Y$ and $X$ are independent, then $ \Omega_{1}=1$, producing the same null distribution used in classical correlation and auto-correlation studies. Compared with \cite{LLT2015}, our result is free of the moment conditions on $X$, while \cite{LLT2015} requires the existence of a fourth order moment on $X$ to achieve the convergence in law.

In order to make statistical inference, such as constructing a confidence interval for $\varrho_{Y, X}(\tau, \iota)$ and testing a hypothesis like $H_0:\varrho_{Y, X}(\tau, \iota)=0$, we need to estimate the covariance function $\Omega_1(\tau_1,\iota_1;\tau_2,\iota_2)$. To this end, denote $m_1(y)=\mathrm{E}\{I(X\leq F_{X}^{-1}(\iota))|Y=y \}$ and $m_2(x)=\mathrm{E}\{ I(Y\leq F_{Y}^{-1}(\tau))|X=x \}$. We can use the nonparametric approach like the Nadaraya-Watson (NW) method (\cite{N1964} and \cite{W1964}) to obtain estimates $\widehat{m}_1(y)$ and $\widehat{m}_2(x)$ for $m_1(y)$ and $m_2(x)$, respectively, where the unknown $F_{X}^{-1}(\iota)$ and $F_{Y}^{-1}(\tau)$ are replaced by $F_{n,X}^{-1}(\iota)$ and $F_{n,Y}^{-1}(\tau)$, respectively.  Therefore, we obtain the estimates $\widehat{\sigma}_{X|Y}(\tau,\iota) = \widehat{m}_1(F_{n,Y}^{-1}(\tau))$ and $\widehat{\sigma}_{Y|X}(\tau,\iota)=\widehat{m}_2(F_{n, X}^{-1}(\iota))$. Next, we give an estimate of $\Omega_{1}(\tau_1,\iota_1;\tau_2,\iota_2)$.
Denote $\widehat{\xi}_n(Y_i,X_i;\tau,\iota)=\big[  I(Y_i \leq F_{n,Y}^{-1}(\tau), X_i\leq F_{n,X}^{-1}(\iota))  - \widehat{\sigma}_{X|Y}(\tau,\iota) I(Y_i\leq F_{n,Y}^{-1}(\tau))  -   \widehat{\sigma}_{Y|X}(\tau,\iota)  I(X_i\leq F_{n,X}^{-1}(\tau))  \big]/\sqrt{\tau(1-\tau) \iota(1-\iota)}$ and $\overline{\xi}_n(Y,X; \tau,\iota)=n^{-1}\sum_{i=1}^{n}\widehat{\xi}_n(Y_i,X_i;\tau,\iota)$. Then, we obtain a consistent estimate of $\Omega_{1}(\tau_1,\iota_1;\tau_2,\iota_2)$ as $\widehat{\Omega}_{1}(\tau_1,\iota_1;\tau_2,\iota_2)= n^{-1}\sum_{i=1}^{n}[\widehat{\xi}_n(Y_i,X_i;\tau_1,\iota_1)-\overline{\xi}_n(Y,X; \tau_1,\iota_1)]\times [\widehat{\xi}_n(Y_i,X_i;\tau_2,\iota_2)-\overline{\xi}_n(Y,X; \tau_2,\iota_2)]$.

In practice, we usually encounter the situation where $Y$ is univariate but $X$ is multivariate. As as extension to Theorem \ref{th1} and to compare the dependence strength of two random variables $X_1$ and $X_2$ on $Y$, we may check the difference  $ \varrho_{Y,X_1}(\tau,\iota)- \varrho_{Y,X_2}(\tau,\iota)$. In particular, we may test a hypothesis by this difference. Given a sample $\{(Y_i,X_{i1}, X_{i2}),i=1,\ldots, n\}$, similarly to (\ref{eq:202b}), we can define $\widehat{\varrho}_{Y, X_1}(\tau, \iota)$ and $\widehat{\varrho}_{Y, X_2}(\tau, \iota)$. The following theorem can be applied to answer this question.

\begin{theorem}  \label{th1b}
Let $0<a<b<1$  and suppose that marginal distributions $F_{X_k}$ and $F_{Y}$ are continuously differentiable on the intervals $[F_{X_k}^{-1}(a)-\varepsilon, F_{X_k}^{-1}(b)+\varepsilon]$ and $[F_{Y}^{-1}(a)-\varepsilon, F_{Y}^{-1}(b)+\varepsilon]$ with positive derivatives $f_{X_k}$ and $f_{Y}$, respectively, for some $\varepsilon >0$ and $k=1,2$. Furthermore, assume that conditional density functions $f_{Y|X_k}$ and $f_{X_k|Y}$, $k=1,2$, are continuous on the product of these intervals. Then we have
\begin{eqnarray*}
   \sqrt{n}\{ [\widehat{\varrho}_{Y, X_1}(\tau, \iota) - \widehat{\varrho}_{Y, X_2}(\tau, \iota)]  -[ \varrho_{Y, X_1}(\tau, \iota)-\varrho_{Y, X_2}(\tau, \iota)] \} \stackrel{w}{\rightsquigarrow } \mathbb{G}_{Y,X_1,X_2}(\tau,\iota),
\end{eqnarray*}
in $\ell^{\infty}([a,b]^2)$, where $\mathbb{G}_{Y,X_1,X_2}(\tau,\iota)$  is Gaussian process with mean zero and covariance function
$\Xi_1(\tau_1,\iota_1;\tau_2,\iota_2) \equiv \mathrm{E}\{ [\eta(Y,X_1,X_2;\tau_1, \iota_1)-\mathrm{E}\eta(Y,X_1,X_2;\tau_1, \iota_1)]\times [\eta(Y,X_1,X_2;\tau_2, \iota_2)-\mathrm{E}\eta(Y,X_1,X_2;\tau_2, \iota_2)]  \}$, $\eta(Y,X_1,X_2;\tau, \iota)= \xi(Y, X_1; \tau, \iota) -  \xi(Y, X_2; \tau, \iota)$ and $\xi(Y, X; \tau, \iota)$ is given in Theorem \ref{th1}.
\end{theorem}

It follows from Theorem \ref{th1b} that for a fixed pair $(\tau, \iota)$, if $\varrho_{Y, X_1}(\tau, \iota)=\varrho_{Y, X_2}(\tau, \iota)$, then $\sqrt{n}\{ [\widehat{\varrho}_{Y, X_1}(\tau, \iota) - \widehat{\varrho}_{Y, X_2}(\tau, \iota)] \stackrel{d}{\to } N(0, \Xi_1)$, where $\Xi_1\equiv \Xi_1(\tau,\iota;\tau,\iota)=\Omega_1^{(1)} + \Omega_1^{(2)} -2A_{12}$, where $\Omega_1^{(k)}$ is the same as $\Omega_1$ except that $X$ involved in $\Omega_1$ is substituted by $X_k$ for $k=1,2$, and
$A_{12} \equiv A_{12}(\tau,\iota)=\big\{ [\sigma_{Y,X_1,X_2}(\tau,\iota) - \sigma_{Y,X_1} \sigma_{Y,X_2}]
- (1-\tau)\sigma_{X_2|Y}\sigma_{Y,X_1}
- \sigma_{Y|X_2}[\sigma_{Y,X_1,X_2}(\tau,\iota)-\iota\sigma_{Y,X_1} ]
- (1-\tau)\sigma_{X_1|Y}\sigma_{Y,X_2}
+ \tau(1-\tau) \sigma_{X_1|Y}\sigma_{X_2|Y}
+ \sigma_{X_1|Y} \sigma_{Y|X_2}(\sigma_{Y,X_2}-\tau\iota)
- \sigma_{Y|X_1}[\sigma_{Y,X_1,X_2}(\tau,\iota)-\iota\sigma_{Y,X_2}]
+ \sigma_{Y|X_1}\sigma_{X_2|Y}(\sigma_{Y,X_1} -\tau\iota)
+ \sigma_{Y|X_1}\sigma_{Y|X_2}[\sigma_{X_1,X_2}(\iota,\iota) -\iota^2]  \big\}/\sqrt{\tau(1-\tau) \iota(1-\iota)}$,
where $\sigma_{Y,X_1,X_2}(\tau,\iota)= F_{Y,X_1,X_2}(F_{Y}^{-1}(\tau), F_{X_1}^{-1}(\iota), F_{X_2}^{-1}(\iota))$ and $\sigma_{X_1,X_2}(\iota,\iota) = F_{X_1,X_2}(F_{X_1}^{-1}(\iota), F_{X_2}^{-1}(\iota))$.
If $Y, X_1$ and $X_2$ are mutually independent, then $\Xi_1=2$. Next, we estimate the covariance function $\Xi_1(\tau_1,\iota_1;\tau_2,\iota_2) $.
Let $\widehat{\eta}_n(Y_i,X_{i1},X_{i2};\tau, \iota)=\widehat{\xi}_n(Y_i, X_{i1}; \tau, \iota) -  \widehat{\xi}_n(Y_i, X_{i2}; \tau, \iota)$, where  $\widehat{\xi}_n(Y_i, X_{i}; \tau, \iota)$ is given before,
and $\overline{\eta}_n(Y,X_{1},X_{2};\tau, \iota)=n^{-1}\sum_{i=1}^{n} \widehat{\eta}_n(Y_i,X_{i1},X_{i2};\tau, \iota)$. Then, $\Xi_1(\tau_1,\iota_1;\tau_2,\iota_2)$ can be estimated as
$\widehat{\Xi}_1(\tau_1,\iota_1;\tau_2,\iota_2)=n^{-1}\sum_{i=1}^{n} [\widehat{\eta}_n(Y_i,X_{i1},X_{i2};\tau_1, \iota_1)- \overline{\eta}_n(Y,X_{1},X_{2};\tau_1, \iota_1)]\times [\widehat{\eta}_n(Y_i,X_{i1},X_{i2};\tau_2, \iota_2)-\overline{\eta}_n(Y,X_{1},X_{2};\tau_2, \iota_2)] $.

\subsection{CC-based Variable Screening}

Suppose that we collect a sample $\{(Y_i, \mathbf{X}_i),i=1,\cdots, n\}$ consisting of $n$ independent copies of $(Y, \mathbf{X})$, where $Y$ is the response variable and $\mathbf{X}=(X_1,\cdots, X_{p})^T$  is a vector of $p$ predictors. When the number, $p$, of predictors is of an exponential order of sample size $n$, i.e., the so-called ultrahigh dimension, and most of $p$ predictors are irrelevant, we can use CC as a screener to identify the sparse set of informative predictors. We write $p_n$ instead of $p$ to emphasize the dependence on sample size. An empirical estimate for CC between $Y$ and $X_j$ is given by
\begin{eqnarray}  \label{eq:202bb}
 \widehat{ \varrho}_{Y, X_j}(\tau, \iota) = \frac{n^{-1}\sum_{i=1}^{n}\psi_{\tau}(Y_i-F_{n, Y}^{-1}(\tau))\psi_{\iota}(X_{ij}-F_{n, X_j}^{-1}(\iota))  }{ \sqrt{\tau(1-\tau)\iota(1-\iota)}}.
\end{eqnarray}
Then, we may select an empirical active set to be
\begin{eqnarray}
  \widehat{\mathcal{M}}_{a} =\big\{ j:  | \widehat{ \varrho}_{Y, X_j}(\tau, \iota) |\geq \nu_n, 1\leq j \leq p_n \big\},  \label{eq:202c}
\end{eqnarray}
 where $\nu_n$ is a user-specified threshold parameter that controls the size of finally screened model. Using CC can lead to sure independence screening (SIS) property and this procedure will be abbreviated as CC-SIS.

Denote the true active set by $\mathcal{M}_a^{\ast}=\{j: | { \varrho}_{Y, X_j}(\tau, \iota) |>0 ,j=1,\ldots, p_n\}$. Write $F_{Y|\mathbf{X}}^{-1}(\tau)=\inf\{y: P(Y\leq y|\mathbf{X})\geq \tau\}$, $u_j= |\varrho_{Y, X_j}(\tau, \iota)|$ and $\widehat{u}_j=|\widehat{ \varrho}_{Y, X_j}(\tau, \iota)|$. To establish the screening consistency, we need the following conditions.
\begin{enumerate}
  \item[(C1)]  In a neighbourhood of $F_{Y}^{-1}(\tau)$, the density $f_{Y}(y)$ of $Y$ is uniformly bounded away from zero and infinity and has bounded derivative. For every $1\leq j\leq p_n$, in a neighbourhood of $F_{X_j}^{-1}(\iota)$, the density $f_{X_j}(x)$ of $X_j$ is uniformly bounded away from zero and infinity and has bounded derivative.
  \item[(C2)] $\min_{j\in\mathcal{M}_a^{\ast}} u_j\geq  C_0 n^{-\kappa}$ for some $\kappa>0$ and $C_0>0$.
\end{enumerate}

\begin{theorem}  \label{th2}
(Screening Property for CC-SIS) Suppose that the condition (C1) holds, \\
(i) for any constant $C>0$, then there exists some positive constant $\tilde{c}_{1}$ such that for sufficiently large $n$,
\begin{eqnarray*}
  P\Big(\max_{1\leq j\leq p_n} \big| \widehat{u}_j - u_j \big| \geq C  n^{-\kappa} \Big) \leq
 6p_n\exp(-\tilde{c}_1n^{1-2\kappa}).
\end{eqnarray*}
(ii) In addition, if condition (C2) is further satisfied and by choosing $\nu_n=C_{1} n^{-\kappa}$ with $C_1 \leq C_0/2$, we have
\begin{eqnarray*}
 &&  P\Big( \mathcal{M}_a^{\ast} \subset  \widehat{ \mathcal{M}}_{a} \Big)  \geq 1- 6s_n \exp(-\tilde{c}_1n^{1-2\kappa})
\end{eqnarray*}
for sufficiently large $n$, where $s_n=|\mathcal{M}_a^{\ast}|$ is the cardinality of set $\mathcal{M}_a^{\ast}$.
\end{theorem}
This result implies that the CC-SIS can select all the truly active predictors with an overwhelming probability. The dimensionality can be as high as $p_n = o(\exp(n^{1-2\kappa}))$, similar to other model-free feature screening methods (see \cite{LPZZ2012} and \cite{WY2015} for example). Moreover, our result requires less condition on both the predictors and the response due to the nonparametric nature. In fact, no moment assumption on the predictors or the response is imposed.

In practice, the threshold parameter $\nu_n$ plays an important role in producing a satisfied model. Small $\nu_n$ value will result in a large number of predictors after screening, which in turn leads to many incorrect positives. Here we consider a data-driven procedure to determine the threshold for the CC-SIS by controlling the false discovery rates (FDR).  By Theorem \ref{th1}, for covariate $j$ such that $ \varrho_{Y, X_j}(\tau, \iota)=0$, it follows that asymptotically, $\sqrt{n} [\widehat{\Omega}_1(\tau,\iota;\tau,\iota)]^{-1/2} \widehat{\varrho}_{Y, X_j}(\tau, \iota) \sim N(0, 1)$. We can use high-criticism $t$-tests to select variables  $\widehat{\mathcal{M}}_{a,\delta}=\{j: \sqrt{n} [\widehat{\Omega}_1(\tau,\iota;\tau,\iota)]^{-1/2} |\widehat{\varrho}_{Y, X_j}(\tau, \iota)| \geq \delta \}$ for a small $\delta>0$. This controls the FDR $\mathrm{E}\{ |\widehat{\mathcal{M}}_{a,\delta} \cap (\mathcal{M}_{a}^{\ast})^{c} |/ | (\mathcal{M}_{a}^{\ast})^{c} | \}$ defined by \cite{ZL2012}. The following proposition justifies this FDR procedure. 


\begin{proposition}  \label{prop1}
(FDR Property) Under conditions (C1)-(C2) and the condition of Theorem \ref{th1}, if we choose $\delta=\Phi^{-1}(1-\bar{d}_n/(2p_n))$ and $\Phi(\cdot)$ is CDF of standard normal variable and $\bar{d}_n$ is the number of false positives that can be tolerated, then for some constant $c_{a}>0$, we have
\begin{eqnarray*}
  \mathrm{E}\Big\{ \frac{|\widehat{\mathcal{M}}_{a, \delta} \cap (\mathcal{M}_{a}^{\ast})^{c} |}{| (\mathcal{M}_{a}^{\ast})^{c} |}   \Big\} \leq \frac{\bar{d}_n}{p_n} + c_a/\sqrt{n}.
\end{eqnarray*}
\end{proposition}

\section{Copula-based Partial Correlation and Variable Screening}
\subsection{Copula-based Partial Correlation, CPC}

To facilitate a joint screening procedure (\cite{MLT2017}), we define a copula-based partial correlation (CPC) for $Y$ and $X$ conditional on a $q$-dimensional random vector $\mathbf{Z}$ as
\begin{eqnarray}
  \varrho_{Y, X|\mathbf{Z} }(\tau, \iota) 
   =    \frac{\mathrm{E}\{\psi_{\tau}(Y-\mathbf{Z}^T\boldsymbol{\alpha}^{0}) \psi_{\iota}(X-\mathbf{Z}^T\boldsymbol{\theta}^{0})\} }{\sqrt{\tau(1-\tau)\iota(1-\iota)}}, \label{eq:203}
\end{eqnarray}
where $\boldsymbol{\alpha}^{0} = \mathrm{argmin}_{\boldsymbol{\alpha}} \mathrm{E}\{\rho_{\tau}(Y-\mathbf{Z}^T\boldsymbol{\alpha})\}$ and $\boldsymbol{\theta}^{0} = \mathrm{argmin}_{\boldsymbol{\theta}} \mathrm{E}\{ \rho_{\iota}(X-\mathbf{Z}^T\boldsymbol{\theta})\}$ and  $\rho_{w}(u)=u[w-I(u\leq 0)]$ for $w=\tau$ or $\iota$. Note that this implies that $\mathbf{Z}^T\boldsymbol{\alpha}^{0} = F_{Y|\mathbf{Z}}^{-1}(\tau)$ and $\mathbf{Z}^T\boldsymbol{\theta}^{0} = F_{X|\mathbf{Z}}^{-1}(\iota)$. Parameters $\ve\alpha$ and $\ve\theta$ can be interpreted as the marginal increment on conditional quantiles of $Y$ and $X$ given $\mathbf{Z}$, respectively, when increasing by a unit of $\mathbf{Z}$. CPC is actually the CC between $Y$ and $X_j$ after removing the confounding effects of $\mathbf{Z}$.  Linear partial correlation has been widely used in regression diagnostics and describes the association of the response and predictor conditional on specifical values of other predictors. The unconditional $\varrho_{Y,X}(\tau,\iota)$ value may be spurious due to lurking variables and does not necessarily imply the same $\varrho_{Y, X|\mathbf{Z}}(\tau, \iota)$ value conditional on $\mathbf{Z}$. Our copula based version is relatively more robust for real data analysis. The CC is a special case of CPC when ${\bf Z}$ is a constant.

With sample observations $\{(Y_i, X_i,\mathbf{Z}_i), i=1,\ldots,n\}$, we can obtain the following estimate of $\varrho_{Y, X|\mathbf{Z}}(\tau, \iota)$. Let
$\widehat{\boldsymbol{\alpha}} = \mathrm{argmin}_{\boldsymbol{\alpha}} \frac{1}{n} \sum_{i=1}^{n} \rho_{\tau}(Y_i-\mathbf{Z}_i^T\boldsymbol{\alpha}) $
and
$\widehat{\boldsymbol{\theta}}= \mathrm{argmin}_{\boldsymbol{\theta}} \frac{1}{n} \sum_{i=1}^{n} \rho_{\iota}(X_{i}-\mathbf{Z}_i^T\boldsymbol{\theta})$. Both can be obtained from a quantile regression straightforwardly. An empirical estimator for $\varrho_{Y, X|\mathbf{Z}}(\tau, \iota)$ is
\begin{eqnarray}
  \widehat{\varrho}_{Y, X|\mathbf{Z}}(\tau, \iota)
   =   \frac{n^{-1}\sum_{i=1}^{n} \psi_{\tau}(Y_i-\mathbf{Z}_i^T\widehat{\boldsymbol{\alpha}}) \psi_{\iota}(X_{i}-\mathbf{Z}_i^T\widehat{\boldsymbol{\theta}}) }{\sqrt{\tau(1-\tau)\iota(1-\iota)}}.   \label{eq:206}
\end{eqnarray}

To study the asymptotic property of $ \widehat{\varrho}_{Y, X|\mathbf{Z}}(\tau, \iota)$, we denote
\begin{eqnarray*}
&&  \Delta_{11} = \mathrm{E}\{ f_{Y|\mathbf{Z}}(\mathbf{Z}^{T}\boldsymbol{\alpha}^{0}) \mathbf{Z}\mathbf{Z}^T \}, \quad
 \Delta_{12} = \mathrm{E}\big\{ F_{X| \mathbf{Z}, Y=\mathbf{Z}^{T}\boldsymbol{\alpha}^{0}}(\mathbf{Z}^{T}\boldsymbol{\theta}^{0}) f_{Y|\mathbf{Z}}(\mathbf{Z}^{T}\boldsymbol{\alpha}^{0}) \mathbf{Z} \big\},  \\
&& \Delta_{21}= \mathrm{E}\big\{ F_{Y| \mathbf{Z}, X=\mathbf{Z}^{T}\boldsymbol{\theta}^{0}}(\mathbf{Z}^{T}\boldsymbol{\alpha}^{0}) f_{X|\mathbf{Z}}(\mathbf{Z}^{T}\boldsymbol{\theta}^{0}) \mathbf{Z}\big\}, \quad
 \Delta_{22} = \mathrm{E}\{ f_{X|\mathbf{Z}}(\mathbf{Z}^{T}\boldsymbol{\theta}^{0}) \mathbf{Z}\mathbf{Z}^T \},     \\
&& \Sigma_{11} = \mathrm{E}\{  F_{Y,X|\mathbf{Z}}(\mathbf{Z}^{T}\boldsymbol{\alpha}^{0}, \mathbf{Z}^{T}\boldsymbol{\theta}^{0})\}[1-\mathrm{E}\{  F_{Y,X|\mathbf{Z}}(\mathbf{Z}^{T}\boldsymbol{\alpha}^{0}, \mathbf{Z}^{T}\boldsymbol{\theta}^{0})\}],   \\
&& \Sigma_{22} = \mathrm{E}\{  \psi_{\tau}^2(Y-\mathbf{Z}^{T}\boldsymbol{\alpha}^{0}) \mathbf{Z}\mathbf{Z}^T\}, \quad
 \Sigma_{33} = \mathrm{E}\{  \psi_{\iota}^2(X-\mathbf{Z}^{T}\boldsymbol{\theta}^{0}) \mathbf{Z}\mathbf{Z}^T\}, \\
&& \Sigma_{12} = \mathrm{E}\{  F_{Y,X|\mathbf{Z}}(\mathbf{Z}^{T}\boldsymbol{\alpha}^{0}, \mathbf{Z}^{T}\boldsymbol{\theta}^{0}) \mathbf{Z} \}, \quad
 \Sigma_{23} = \mathrm{E}\{  \psi_{\tau}(Y-\mathbf{Z}^{T}\boldsymbol{\alpha}^{0})\psi_{\iota}(X-\mathbf{Z}^{T}\boldsymbol{\theta}^{0}) \mathbf{Z}\mathbf{Z}^T\},
\end{eqnarray*}
where $\boldsymbol{\alpha}^{0}$ and $\boldsymbol{\theta}^{0}$ are defined in (\ref{eq:203}).

We have the following asymptotic result for the CPC.

\begin{theorem}  \label{th3}
Let $0<a<b<1$. Suppose that $\Delta_{11}$ and $\Delta_{22}$ are uniformly positive definite matrices in $\tau$ and $\iota$, and there exists a constant $\pi>0$ such that $f_{Y|\mathbf{Z}}(\mathbf{Z}^{T} \boldsymbol{\alpha}^{0}+\cdot)$, $f_{Y|\mathbf{Z}, X}(\mathbf{Z}^{T} \boldsymbol{\alpha}^{0}+\cdot) $, $f_{X|\mathbf{Z}, Y}(\mathbf{Z}^{T} \boldsymbol{\theta}^{0}+\cdot)$ and  $f_{X|\mathbf{Z}}(\mathbf{Z}^{T} \boldsymbol{\theta}^{0}+\cdot)$ are uniformly integrable on $[-\pi, \pi]$ and uniformly bounded away from zero and infinity in $\tau$ and $\iota$. Then
\begin{eqnarray*}
    \sqrt{n}\{ \widehat{\varrho}_{Y, X|\mathbf{Z}}(\tau, \iota) -   \varrho_{Y, X|\mathbf{Z}}(\tau, \iota ) \} \stackrel{w}{\rightsquigarrow} \mathbb{G}_{Y,X|\mathbf{Z}}(\tau, \iota)
\end{eqnarray*}
in $\ell^{\infty}([a,b]^2)$, where $\mathbb{G}_{Y,X|\mathbf{Z}}(\tau,\iota)$ is Gaussian process with mean zero and covariance function
$\Omega_{2}(\tau_1,\iota_1;\tau_2,\iota_2) \equiv \mathrm{E}\{[\zeta (Y,X, \mathbf{Z}; \tau_1, \iota_1) - \mathrm{E} \zeta(Y, X, \mathbf{Z}; \tau_1, \iota_1)]\times [\zeta(Y, X, \mathbf{Z}; \tau_2, \iota_2) - \mathrm{E} \zeta(Y,X, \mathbf{Z}; \tau_2, \iota_2)]\}$ and $\zeta (Y,X, \mathbf{Z}; \tau, \iota)  =  [ I(Y\leq \mathbf{Z}^{T}\boldsymbol{\alpha}^{0}, X\leq  \mathbf{Z}^{T}\boldsymbol{\theta}^{0} ) -  \Delta_{12}^T\Delta_{11}^{-1} I(Y\leq \mathbf{Z}^{T}\boldsymbol{\alpha}^{0}) \mathbf{Z} - \Delta_{21}^T\Delta_{22}^{-1} I(X\leq \mathbf{Z}^{T}\boldsymbol{\theta}^{0}) \mathbf{Z}]/\sqrt{ \tau(1-\tau) \iota(1-\iota)} $.
\end{theorem}

If $\mathbf{Z}\equiv 1$, in another word, there is no conditional variable available, the asymptotic distribution in Theorem \ref{th3} reduces to that in Theorem \ref{th1}. The above result implies that for a fixed pair $(\tau,\iota)$, if $\varrho_{Y, X|\mathbf{Z}}(\tau, \iota )=0$, then $\sqrt{n} \widehat{\varrho}_{Y, X|\mathbf{Z}}(\tau, \iota) \stackrel{d}{\to}  N(0, \Omega_2)$, where  $ \Omega_2  \equiv \Omega_2(\tau,\iota;\tau,\iota) = \mathrm{E}\{[\zeta (Y,X, \mathbf{Z}; \tau, \iota) - \mathrm{E} \zeta(Y, X, \mathbf{Z}; \tau, \iota)]^2\} = \frac{1}{\tau(1-\tau) \iota(1-\iota)} [   \Sigma_{11} + \Delta_{12}^T \Delta_{11}^{-1} \Sigma_{22} \Delta_{11}^{-1} \Delta_{12}
+ \Delta_{21}^T \Delta_{22}^{-1} \Sigma_{33} \Delta_{22}^{-1} \Delta_{21}
- 2(1-\tau) \Delta_{12}^T \Delta_{11}^{-1}\Sigma_{12} - 2(1-\iota) \Delta_{21}^T \Delta_{22}^{-1} \Sigma_{12} + 2\Delta_{12}^T \Delta_{11}^{-1} \Sigma_{23}\Delta_{22}^{-1}\Delta_{21} ]$.

This theorem can be used for statistical inference if we can find a consistent estimate of $\Omega_2$. To this end, let $e_1^{\ast}=Y-\mathbf{Z}^{T}\boldsymbol{\alpha}^{0}$ and $e_2^{\ast}=X-\mathbf{Z}^{T}\boldsymbol{\theta}^{0}$ and assume that the random vectors $(e_1^{\ast},\mathbf{Z},X)$ and $(e_2^{\ast},\mathbf{Z},Y)$ have joint densities $f_{e_1^{\ast},\mathbf{Z},X}$ and $f_{e_2^{\ast},\mathbf{Z},Y}$, respectively. Denote by $f_{e_1^{\ast}}$, $f_{e_2^{\ast}}$, $f_{e_1^{\ast}|\mathbf{Z}}$, $f_{e_1^{\ast}|\mathbf{Z}, X}$, $f_{e_2^{\ast}|\mathbf{Z}}$ and $f_{e_2^{\ast}|\mathbf{Z}, Y}$ the marginal densities of $e_1^{\ast}$ and of $e_2^{\ast}$, the conditional densities of $e_1^{\ast}$ given $\mathbf{Z}$ and $(\mathbf{Z},X)$ and of $e_2^{\ast}$ given $\mathbf{Z}$ and $(\mathbf{Z}, Y)$, respectively. Then, it can be verified that
$\Delta_{11}=\mathrm{E}\{ f_{e_1^{\ast}|\mathbf{Z}}(0) \mathbf{Z}\mathbf{Z}^T \} = f_{e_1^{\ast}}(0)\mathrm{E}\{ \mathbf{Z}\mathbf{Z}^T  | e_1^{\ast}=0\}$
and, similarly,
$\Delta_{12} =f_{e_1^{\ast}}(0)\mathrm{E}\{ I(X\leq \mathbf{Z}^{T}\boldsymbol{\theta}^{0}) \mathbf{Z}  | e_1^{\ast}=0\}  $,
$\Delta_{21}=f_{e_2^{\ast}}(0) \mathrm{E}\{ I(Y\leq \mathbf{Z}^{T}\boldsymbol{\alpha}^{0}) \mathbf{Z}| e_2^{\ast}=0\}$
and
$\Delta_{22}=f_{e_2^{\ast}}(0) \mathrm{E}\{ \mathbf{Z}\mathbf{Z}^T |e_2^{\ast}=0 \}$. To estimate these quantities, we first calculate the quantile regression estimates $\widehat{\ve\alpha}$ and $\widehat{\ve\theta}$ and then obtain the corresponding quantile residuals
$\hat{e}_{1i}^{\ast}=Y_i-\mathbf{Z}_i^T\widehat{\ve\alpha}$ and $\hat{e}_{2i}^{\ast}=X_i-\mathbf{Z}_i^T\widehat{\ve\theta}$ for $i=1,\ldots, n$.
Next, we provide an estimate for $\Delta_{12}$ and the estimators for $\Delta_{11}$, $\Delta_{21}$ and $\Delta_{22}$ can be obtained similarly.
We can use nonparametric NW estimates used in estimating $\sigma_{X|Y}(\tau,\iota)$ and $\sigma_{Y|X}(\tau,\iota)$ in Section 2.1 to obtain estimates for each component of $\mathbf{m}(s)=\mathrm{E}\{ I(X\leq \mathbf{Z}^{T}\boldsymbol{\theta}^{0}) \mathbf{Z}  | e_1^{\ast}=s\} $ using the data $\{ (\hat{e}_{1i}^{\ast}, I(X_i\leq \mathbf{Z}_i^T\widehat{\boldsymbol{\theta}}) \mathbf{Z}_i),i=1,\ldots, n\}$, and denote it by $\widehat{\mathbf{m}}(s)$. Then we obtain  $\widehat{\Delta}_{12}=\widehat{f}_{e_1^{\ast}}(0)\widehat{\mathbf{m}}(0) $,  where $\widehat{f}_{e_1^{\ast}}(0)$ is a nonparametric kernel density estimate for $f_{e_1^{\ast}}(0)$ in $\Delta_{12}$ based on $\{\hat{e}_{1i}^{\ast}, i=1,\ldots, n\}$. It can be shown that such $\widehat{\Delta}_{12}$ is consistent under some regularity conditions.
For other unknown terms involved in $\Omega_2$, we have
$\widehat{\Sigma}_{11} = n^{-1}\sum_{i=1}^{n}  I(Y_i\leq \mathbf{Z}_i^T\widehat{\boldsymbol{\alpha}}, X_i\leq \mathbf{Z}_i^T\widehat{\boldsymbol{\theta}})-[ n^{-1}\sum_{i=1}^{n}  I(Y_i\leq \mathbf{Z}_i^T\widehat{\boldsymbol{\alpha}}, X_i\leq \mathbf{Z}_i^T\widehat{\boldsymbol{\theta}})]^2$,
$\widehat{\Sigma}_{22} = n^{-1}\sum_{i=1}^{n} \psi_{\tau}^2(Y_i-\mathbf{Z}_i^T\widehat{\boldsymbol{\alpha}}) \mathbf{Z}_i\mathbf{Z}_i^T $,
$\widehat{\Sigma}_{33} = n^{-1}\sum_{i=1}^{n}   \psi_{\iota}^2(X_i-\mathbf{Z}_i^T\widehat{\boldsymbol{\theta}}) \mathbf{Z}_i\mathbf{Z}_i^T$,
$\widehat{\Sigma}_{12} = n^{-1}\sum_{i=1}^{n} I(Y_i\leq \mathbf{Z}_i^T\widehat{\boldsymbol{\alpha}}, X_i\leq \mathbf{Z}_i^T\widehat{\boldsymbol{\theta}}) \mathbf{Z}_i $,
and
$\widehat{\Sigma}_{23} = n^{-1}\sum_{i=1}^{n} \psi_{\tau}(Y_i-\mathbf{Z}_i^T\widehat{\boldsymbol{\alpha}})\psi_{\iota}(X_i-\mathbf{Z}_i^T\widehat{\boldsymbol{\theta}}) \mathbf{Z}_i\mathbf{Z}_i^T$.
Using the plug-in approach, a consistent estimate of $\Omega_2$ is thus obtained and denoted by $\widehat{\Omega}_2$.

The next theorem can be used to test whether $\varrho_{\tau, \iota}(Y, X_1|\mathbf{Z}) = \varrho_{\tau, \iota}(Y, X_2|\mathbf{Z})$
for two different random variables $X_1$ and $X_2$. Write $\ve\theta_k^{0}=\mathrm{argmin}_{\boldsymbol{\theta}} \mathrm{E}\{ \rho_{\iota} (X_k - \mathbf{Z}^T\ve\theta) \}$ for $k=1,2$ and let $\Delta_{12}^{(k)}$ be $\Delta_{12}$, where the involved $X$ and $\ve\theta^{0}$ are replaced by $X_k$ and $\ve\theta_k^{0}$, respectively, for $k=1,2$. In the same manner, we can define $\Delta_{21}^{(k)}$, $ \Delta_{22}^{(k)}$, $\Sigma_{11}^{(k)}$, $\Sigma_{33}^{(k)}$, $\Sigma_{12}^{(k)}$, $\Sigma_{23}^{(k)}$ and, accordingly, $\Omega_{2}^{(k)}$ with $k=1,2$. In addition, we write
 $ \Delta_{31}=\mathrm{E}\{ F_{Y, X_1, X_2|\mathbf{Z}}(\mathbf{Z}^{T}\boldsymbol{\alpha}^{0}, \mathbf{Z}^T\ve\theta_1^{0}, \mathbf{Z}^T\ve\theta_2^{0}) \}
  - \mathrm{E}\{ F_{Y,X_1|\mathbf{Z}}( \mathbf{Z}^{T}\boldsymbol{\alpha}^{0},  \mathbf{Z}^T\ve\theta_1^{0}) \}\mathrm{E}\{ F_{Y,X_2|\mathbf{Z}} (\mathbf{Z}^{T}\boldsymbol{\alpha}^{0}, \mathbf{Z}^T\ve\theta_2^{0})\}$,
 $ \Delta_{32}=\mathrm{E}\{ F_{Y,X_1, X_2|\mathbf{Z}}(\mathbf{Z}^{T}\boldsymbol{\alpha}^{0}, \mathbf{Z}^T\ve\theta_1^{0}, \mathbf{Z}^T\ve\theta_2^{0})\mathbf{Z} \}$ and
$\Delta_{33} = \mathrm{E}\{ \psi_{\iota}(X_1 - \mathbf{Z}^T\ve\theta_1^{0})  \psi_{\iota}(X_2 - \mathbf{Z}^T\ve\theta_2^{0}) \mathbf{Z}\mathbf{Z}^T\}$.

\begin{theorem}  \label{th3b}
Let $0<a<b<1$. Suppose that matrices $\Delta_{11}$ and $\Delta_{22}^{(k)}, k=1,2$, are uniformly positive definite in $\tau$ and $\iota$, and there exists a constant $\pi>0$ such that $f_{Y|\mathbf{Z}}(\mathbf{Z}^{T} \boldsymbol{\alpha}^{0}+\cdot)$, $f_{Y|\mathbf{Z}, X_k}(\mathbf{Z}^{T} \boldsymbol{\alpha}^{0}+\cdot) $, $f_{X_k|\mathbf{Z}, Y}(\mathbf{Z}^{T} \boldsymbol{\theta}^{0}+\cdot)$ and  $f_{X_k|\mathbf{Z}}(\mathbf{Z}^{T} \boldsymbol{\theta}^{0}+\cdot)$ are uniformly integrable on $[-\pi, \pi]$ for $k=1,2$ and uniformly bounded away from zero and infinity in $\tau$ and $\iota$. Then
\begin{eqnarray*}
\sqrt{n}\{ [\widehat{\varrho}_{Y, X_1|\mathbf{Z}}(\tau, \iota) -\widehat{\varrho}_{Y, X_2|\mathbf{Z}}(\tau, \iota)]-   [\varrho_{Y, X_1|\mathbf{Z}}(\tau, \iota)-\varrho_{Y, X_2|\mathbf{Z}}(\tau, \iota)] \} \stackrel{w}{\rightsquigarrow} \mathbb{G}_{Y,X_1, X_2|\mathbf{Z}}(\tau, \iota)
\end{eqnarray*}
in $\ell^{\infty}([a,b]^2)$, where $\mathbb{G}_{Y,X_1, X_2|\mathbf{Z}}(\tau,\iota)$ is Gaussian process with mean zero and covariance function
$\Xi_{2}(\tau_1,\iota_1; \tau_2,\iota_2) \equiv \mathrm{E}\{[\beta (Y,X_1,X_2, \mathbf{Z}; \tau_1, \iota_1) - \mathrm{E} \beta(Y, X_1,X_2, \mathbf{Z}; \tau_1, \iota_1)]\times [\beta(Y, X_1,X_2, \mathbf{Z}; \tau_2, \iota_2) - \mathrm{E} \beta(Y,X_1, X_2, \mathbf{Z}; \tau_2, \iota_2)]\}$ and $\beta (Y,X_1, X_2, \mathbf{Z}; \tau, \iota)  =  \zeta (Y,X_1, \mathbf{Z}; \tau, \iota) - \zeta (Y,X_2, \mathbf{Z}; \tau, \iota)$, where $\zeta (Y,X_1, \mathbf{Z}; \tau, \iota)$ is given in Theorem \ref{th3}.
\end{theorem}

 For fixed $(\tau, \iota)$, if $\varrho_{\tau, \iota}(Y, X_1|\mathbf{Z}) = \varrho_{\tau, \iota}(Y, X_2|\mathbf{Z})$, then  $\sqrt{n} [\widehat{\varrho}_{Y, X_1|\mathbf{Z}}(\tau, \iota) -\widehat{\varrho}_{Y, X_2|\mathbf{Z}}(\tau, \iota)] \stackrel{d}{\to} N(0, \Xi_2)$, where $\Xi_2\equiv \Xi_{2}(\tau,\iota; \tau,\iota) = \Omega_2^{(1)}+\Omega_2^{(2)} -2B_{12}$ and $ B_{12} \equiv B_{12}(\tau,\iota)= \frac{1}{\tau(1-\tau) \iota(1-\iota)} [ \Delta_{31} - (1-\tau)(\Delta_{12}^{(2)})^T\Delta_{11}^{-1} \Sigma_{12}^{(1)}  -  (1-\tau)(\Delta_{12}^{(1)})^T\Delta_{11}^{-1} \Sigma_{12}^{(2)}
+  (\Delta_{12}^{(1)})^T\Delta_{11}^{-1} \Sigma_{22}\Delta_{11}^{-1} \Delta_{12}^{(2)} + (\Delta_{12}^{(1)})^T  \Delta_{11}^{-1} \times \Sigma_{23}^{(2)} (\Delta_{22}^{(2)})^{-1} \Delta_{21}^{(2)}
+ (\Delta_{21}^{(1)})^T(\Delta_{22}^{(1)})^{-1} \Delta_{32} + (\Delta_{21}^{(1)})^T(\Delta_{22}^{(1)})^{-1} \Sigma_{23}^{(1)} \Delta_{11}^{-1}\Delta_{12}^{(2)}
+ (\Delta_{21}^{(2)})^T(\Delta_{22}^{(2)})^{-1}(\iota \Sigma_{12}^{(1)} -\Delta_{32}) +  (\Delta_{21}^{(1)})^T(\Delta_{22}^{(1)})^{-1}  \Delta_{33} (\Delta_{22}^{(2)})^{-1}\Delta_{21}^{(2)}]$.
Given a sample of observations $\{(Y_i, X_{i1}, X_{i2},\mathbf{Z}_i), i=1,\ldots, n\}$, the asymptotic variance $\Xi_2$ can be estimated as $\widehat{\Xi}_2 = \widehat{\Omega}_2^{(1)} + \widehat{\Omega}_{2}^{(2)} - 2\widehat{B}_{12}$, where $\widehat{\Omega}_2^{(1)}$ and $\widehat{\Omega}_2^{(2)}$ are defined as $\widehat{\Omega}_2$ given above. To obtain the estimate $\widehat{B}_{12}$, we only need to estimate $\Delta_{31}$, $\textcolor[rgb]{1.00,0.00,0.00}{\Delta_{32}}$ and $\textcolor[rgb]{1.00,0.00,0.00}{\Delta_{33}}$ since the rest of unknown quantities involved in $B_{12}$ can be estimated using the previous methods. To this end, we can use the following estimates:
$\widehat{\Delta}_{31} =n^{-1}\sum_{i=1}^{n} I(Y_i\leq \mathbf{Z}_i^T\widehat{\ve\alpha}, X_{i1} \leq \mathbf{Z}_i^T\widehat{\ve\theta}_1, X_{i2} \leq \mathbf{Z}_i^T\widehat{\ve\theta}_2)  - \big[n^{-1}\sum_{i=1}^{n} I(Y_i\leq \mathbf{Z}_i^T\widehat{\ve\alpha}, X_{i1} \leq \mathbf{Z}_i^T\widehat{\ve\theta}_1)\big] \big[n^{-1}\sum_{i=1}^{n} I(Y_i\leq \mathbf{Z}_i^T\widehat{\ve\alpha}, X_{i2} \leq \mathbf{Z}_i^T\widehat{\ve\theta}_2)\big]$,
$\widehat{\Delta}_{32}= n^{-1}\sum_{i=1}^{n} I(Y_i\leq \mathbf{Z}_i^T\widehat{\ve\alpha}, X_{i1} \leq \mathbf{Z}_i^T\widehat{\ve\theta}_1, X_{i2}\leq \mathbf{Z}_i^T\widehat{\ve\theta}_2) \mathbf{Z}_i$,
$\widehat{\Delta}_{33}= n^{-1}\sum_{i=1}^{n} \psi_{\iota}(X_{i1} - \mathbf{Z}_i^T\widehat{\ve\theta}_1)  \psi_{\iota}(X_{i2} - \mathbf{Z}_i^T\widehat{\ve\theta}_2) \mathbf{Z}_i\mathbf{Z}_i^T$,
where
$\widehat{\ve\theta}_k= \mathrm{argmin}_{\boldsymbol{\theta}} n^{-1}\sum_{i=1}^{n} \rho_{\iota}(X_{ik}-\mathbf{Z}_{i}^T \ve\theta), k=1,2$.

\subsection{CPC-based Variable Screening}
We may now propose a joint robust screening using the CPC. There are two practical scenarios to favor joint screening over marginal screening. First, very often we may acquire low-dimensional variables $\mathbf{W} \in \mathbb{R}^{r}$ in addition to ultrahigh dimensional covariates $\mathbf{X}$. For example, when studying the relationship between a disease phenotype $Y$ and genetic variables ${\bf X}$, we may also have patient demographical information or environmental variables and include them in $\mathbf{W}$. Consequently, we have a data set $\{(Y_i, \mathbf{X}_i, \mathbf{W}_i), i=1,\ldots, n\}$. Second, even if there is \textcolor[rgb]{1.00,0.00,0.00}{no external variable} ${\bf W}$, it may still be necessary to consider a joint screening by removing the effects from correlated components in ${\bf X}$. For instance, some covariates, $\mathbf{X}_{\mathcal{S}_j}$, may be closely correlated to $X_j$ and influence the observed correlation between $Y$ and $X_j$ indirectly, where $\mathcal{S}_j$ is a subset of $\{1,\cdots, p_n\}\setminus \{j\}$. \cite{MLT2017} also considered a set $\mathcal{S}_j$ which is referred to as a conditional set with relatively small size ($<n$). To account for both scenarios, we may consider the conditional variables $\mathbf{Z}=(\mathbf{W}^T,\mathbf{X}_{\mathcal{S}_j}^T)^T$ in this paper. We allow that the conditional variables $\mathbf{Z}$ to differ with $j$. However, for simplicity of presentation, we still use $\mathbf{Z}$ instead of $\mathbf{Z}_j$, and we denote by $q_n$ the dimension of $\mathbf{Z}$. In principle, we only need $q_n= \max_{1\leq j\leq p_n} (r+|\mathcal{S}_j|)$ for sure screening. In practice, we may select a proper $\mathcal{S}_j$ as follows: Treat $X_j$ as the response and $\mathbf{X}_{-j}=\{X_k, k\neq j, 1\leq k\leq p_n\}$ as the predictors and then apply any sensible marginal screening method such as the CC-SIS to pick out the top $\ell$ most important predictors and set them as the conditional variables.

 For ultrahigh dimensional covariates $\mathbf{X}=(X_1,\cdots, X_{p_n})^{T}$, we can define CPC between $Y$ and the $j$th marker $X_j$ given $\mathbf{Z}$ in the same way as in (\ref{eq:203}), namely,
\begin{eqnarray}
  \varrho_{Y, X_j|\mathbf{Z}}(\tau, \iota)  =   \frac{\mathrm{E}\{\psi_{\tau}(Y-\mathbf{Z}^{T}\boldsymbol{\alpha}^{0}) \psi_{\iota}(X_j-\mathbf{Z}^T\boldsymbol{\theta}_j^{0})\} }{\sqrt{\tau(1-\tau)\iota(1-\iota)}}, \label{eq:203b}
\end{eqnarray}
where $\boldsymbol{\alpha}^{0} = \mathrm{argmin}_{\boldsymbol{\alpha}} \mathrm{E}\{\rho_{\tau}(Y-\mathbf{Z}^T\boldsymbol{\alpha})\}$ and $\boldsymbol{\theta}_j^{0} = \mathrm{argmin}_{\boldsymbol{\theta}_j} \mathrm{E}\{ \rho_{\iota}(X_j-\mathbf{Z}^T\boldsymbol{\theta}_j)\}$.  As in (\ref{eq:206}),  a sample estimate for $\varrho_{Y, X_j|\mathbf{Z}}(\tau, \iota)$ can be given as
\begin{eqnarray}
  \widehat{\varrho}_{Y, X_j|\mathbf{Z}}(\tau, \iota)
   =   \frac{n^{-1}\sum_{i=1}^{n} \psi_{\tau}(Y_i-\mathbf{Z}_i^T\widehat{\boldsymbol{\alpha}}) \psi_{\iota}(X_{ij}-\mathbf{Z}_i^T\widehat{\boldsymbol{\theta}}_j) }{\sqrt{\tau(1-\tau)\iota(1-\iota)}},   \label{eq:206b}
\end{eqnarray}
where
$\widehat{\boldsymbol{\alpha}} = \mathrm{argmin}_{\boldsymbol{\alpha}} \frac{1}{n} \sum_{i=1}^{n} \rho_{\tau}(Y_i-\mathbf{Z}_i^T\boldsymbol{\alpha})$
and
$\widehat{\boldsymbol{\theta}}_j= \mathrm{argmin}_{\boldsymbol{\theta}_j} \frac{1}{n} \sum_{i=1}^{n} \rho_{\iota}(X_{ij}-\mathbf{Z}_i^T\boldsymbol{\theta}_j)$.
The CPC screening yields the following empirical active set:
\begin{eqnarray}
  \widehat{\mathcal{M}}_{b} =\{ j:  |\widehat{ \varrho}_{Y, X_j|\mathbf{Z} }(\tau, \iota )|\geq v_n, 1\leq j \leq {p_n} \},  \label{eq:207}
\end{eqnarray}
where $v_n$ is a user-specified threshold parameter. We refer to this sure independence screening procedure as CPC-SIS. Clearly, CPC-SIS extends earlier conditional sure independence screening such as \cite{BFV2016}.

Let $\mathcal{M}_b^{\ast}=\{j: |{ \varrho}_{Y, X_j|\mathbf{Z} }(\tau, \iota )|>0 , j=1,\ldots, p\}$ be the true active set. We write $F_{Y|\mathbf{X}, \mathbf{W}}^{-1}(\tau)=\inf\{y: P(Y\leq y|\mathbf{X},\mathbf{W})\geq \tau\}$. For simplicity, we still use $u_j=| \varrho_{Y, X_j|\mathbf{Z}}(\tau, \iota)|$ and $\widehat{u}_j = |\widehat{ \varrho}_{Y, X_j|\mathbf{Z}}(\tau, \iota )|$ to denote the underlying and empirical CPC utilities, respectively.  To establish the sure screening property, we need the following conditions, which are very mild and similarly imposed in \cite{MLT2017}.
\begin{enumerate}
  \item[(D1)] (i) The conditional density $f_{Y|\mathbf{Z}=\mathbf{z}}(y)$ of $Y$ given $\mathbf{Z}=\mathbf{z}$ satisfies the Lipschitz condition of order $1$ and $f_{Y|\mathbf{Z}=\mathbf{z}}(y)>0$ for any $y$ in a neighborhood of $\mathbf{Z}^{T}\boldsymbol{\alpha}^{0}=\mathbf{z}^{T}\boldsymbol{\alpha}^{0}$.  \\
      (ii) For every $1\leq j\leq p_n$, the conditional density $f_{X_j|\mathbf{Z}=\mathbf{z}}(x)$ of $X_j$ given $\mathbf{Z}=\mathbf{z}$ satisfies the Lipschitz condition of order $1$ and $f_{X_j|\mathbf{Z}=\mathbf{z}}(x)>0$ for any $x$ in a neighborhood of $\mathbf{Z}^{T}\boldsymbol{\theta}^{0}=\mathbf{z}^{T}\boldsymbol{\theta}^{0}$.
  \item[(D2)] (i) There exist some finite constants $m_1, m_2$ and $m_3$ such that
  \[  \max_{i,j} |Z_{ij}| \leq m_1, \quad \max_{i} |\mathbf{Z}_i^{T}\ve\alpha^0| \leq m_2, \quad \max_{i,j} |\mathbf{Z}_i^{T}\ve\theta_j^0| \leq m_3.
  \]
  (ii) There exist two positive finite constants $c_{\min}$ and $c_{\max}$ such that
  \[  c_{\min} \leq \lambda_{\min}(\mathrm{E}(\mathbf{Z}\mathbf{Z}^T)) \leq \lambda_{\max}(\mathrm{E}(\mathbf{Z}\mathbf{Z}^T)) \leq c_{\max},
  \]
  where $\lambda_{\min}(\mathrm{E}(\mathbf{Z}\mathbf{Z}^T))$ and $\lambda_{\max}(\mathrm{E}(\mathbf{Z}\mathbf{Z}^T))$ stand for the minimum and maximum eigenvalues of $\mathrm{E}(\mathbf{Z}\mathbf{Z}^T)$, respectively.
  \item[(D3)] $\min_{j\in\mathcal{M}_b^{\ast}} u_j\geq  C_{0}^{\ast} n^{-\kappa}$ for some $\kappa>0$ and $C_0^{\ast}>0$.
\end{enumerate}

\begin{theorem}  \label{th4}
(Screening Property for CPC-SIS) Suppose that the conditions (D1) and (D2) hold, \\
(i) for any constant $C>0$, then there exists some positive constant $\tilde{c}_{1}^{\ast}$ such that for sufficiently large $n$,
\begin{eqnarray*}
  P\Big(\max_{1\leq j\leq p_n} \big| \widehat{u}_j - u_j \big| \geq C  n^{-\kappa} \Big) \leq
 12p_n\exp(-\tilde{c}_1^{\ast} q_n^{-1} n^{1-2\kappa}).
\end{eqnarray*}
(ii) In addition, if condition (D3) is further satisfied and by choosing $v_n=C_{2} n^{-\kappa}$ with $C_2 \leq C_0^{\ast}/2$, we have
\begin{eqnarray*}
 &&  P\Big( \mathcal{M}_b^{\ast} \subset  \widehat{ \mathcal{M}}_{b} \Big)  \geq 1- 12s_n \exp(-\tilde{c}_1^{\ast} q_n^{-1} n^{1-2\kappa})
\end{eqnarray*}
for sufficiently large $n$, where $s_n=|\mathcal{M}_b^{\ast}|$.
\end{theorem}

When conditional variables are available, our proposed CPC-SIS method can handle the dimensionality of order $ p_n = o(\exp(q_n^{-1} n^{1-2\kappa}))$. If $q_n=O(1)$, then the dimension can be as high as $o(n^{1-2\kappa})$, the same order as that of the CC-SIS. Moreover, the proposed CPC-SIS can be readily used for the ultrahigh dimensional data as long as $q_n=o(n^{1-2\kappa})$.

As in Section 2.2, we can determine a proper $v_n$ by controlling FDR. By Theorem \ref{th3}, for covariate $j$ such that $ \varrho_{Y, X_j|\mathbf{Z}}(\tau, \iota)=0$, we have $\sqrt{n} \widehat{\Omega}_2^{-1/2} \widehat{\varrho}_{Y, X_j|\mathbf{Z}}(\tau, \iota) \sim N(0, 1)$ asymptotically. Then, we select variables $\widehat{\mathcal{M}}_{b,\delta}=\{j: \sqrt{n} \widehat{\Omega}_2^{-1/2} |\widehat{\varrho}_{Y, X_j|\mathbf{Z}}(\tau, \iota)| \geq \delta \}$ for a small $\delta>0$, which controls the FDR $\mathrm{E}\{ |\widehat{\mathcal{M}}_{b,\delta} \cap (\mathcal{M}_{b}^{\ast})^{c} |/ | (\mathcal{M}_{b}^{\ast})^{c} | \}$. 

\begin{proposition}  \label{prop2}
(FDR Property) Under conditions (D1)-(D3) and the condition of Theorem \ref{th3}, if we choose $\delta=\Phi^{-1}(1-\bar{d}_n/(2p_n))$ and $\Phi(\cdot)$ and $\bar{d}_n$ are the same as those in Proposition \ref{prop1}, then for some constant $c_{b}>0$, we have
\begin{eqnarray*}
  \mathrm{E}\Big\{ \frac{|\widehat{\mathcal{M}}_{b,\delta} \cap (\mathcal{M}_{b}^{\ast})^{c} |}{| (\mathcal{M}_{b}^{\ast})^{c} | }  \Big\} \leq \frac{\bar{d}_n}{p_n} + c_b/\sqrt{n}.
\end{eqnarray*}
\end{proposition}

\section{Implementation on CPC-SIS}
We provide more details on the implementation of CPC-SIS. We consider three practical types of conditional variables in the following.

\textit{Case 1. } If $\mathbf{W}$ is not available, we consider the conditional variables from $\mathbf{X}$ itself for each $X_j$, namely, $\mathbf{Z}=\mathbf{X}_{\mathcal{S}_j}$ for $j=1,\ldots, p_n$. We start with an empty active set $\mathcal{A}^{(0)}=\emptyset$.
\begin{itemize}
  \item Step 1. For $j=1,\ldots, p_n$, select confounding sets $\mathcal{S}_j^{\nu}$'s via the partial correlation based consequential test (\cite{MLT2017}).
  \item Step 2. In the $k$th iteration, where $k=1, \ldots, d^{*}$ and $d^{*}=\lfloor 2 (n/\log n)^{1/2} \rfloor$, for given $\mathcal{A}^{(k-1)}$, we update $\mathcal{S}_j=\mathcal{A}^{(k-1)} \cup \mathcal{S}_j^{\nu}$ and then find the variable index  $j^{*}$ such that $j^{*} = \mathrm{argmax}_{j\not\in \mathcal{A}^{(k-1)} }| \widehat{\varrho}_{Y, X_j|\mathbf{Z}}(\tau, \iota)|$. Update $\mathcal{A}^{(k)}=\mathcal{A}^{(k-1)}\cup \{j^{*}\}$.
  \item Step 3. In the $k$th iteration, where  $k=d^{*}+1, \ldots, d_n$,  we set $\mathcal{S}_j=\mathcal{A}^{(d^{*})} \cup \mathcal{S}_j^{\nu}$ and then find $j^{*} = \mathrm{argmax}_{j\not\in \mathcal{A}^{(k-1)} }| \widehat{\varrho}_{Y, X_j|\mathbf{Z}}(\tau, \iota)|$. Update $\mathcal{A}^{(k)}=\mathcal{A}^{(k-1)}\cup \{j^{*}\}$. Use $\mathcal{A}^{(d_n)}\equiv \widehat{\cal M}_b$ as the final set of selected covariates.
\end{itemize}
It is worth noting that the main difference between Steps 2 and 3 is that for Step 2, the conditional set is updated gradually via adding one selected index variable in the first $d^{\ast}$ iterations, while for Step 3 the conditional set keeps intact in the last $d_n-d^{\ast}$ iterations.

\textit{Case 2. } If $\mathbf{W}$ is available, we consider the same conditional variables for each target $X_j$, namely, $\mathbf{Z}=\mathbf{W}$ for $j=1,\ldots, p_n$.
\begin{itemize}
  \item Step 1. For $j=1,\ldots, p_n$, compute the CPC utility statistics $\widehat{u}_j = |\widehat{\varrho}_{Y, X_j|\mathbf{Z}}(\tau, \iota)|$.
  \item Step 2. Rank the covariates in terms of their $\widehat{u}_j$'s in a decreasing order and then select the top $d_n$ covariates as the final set of selected covariates.
\end{itemize}

\textit{Case 3. } If $\mathbf{W}$ is available, we slightly modify the algorithm in Case 1. The steps are the same as those in Case 1 only except that we consider the conditional variables $\mathbf{Z}=(\mathbf{W}^T, \mathbf{X}_{\mathcal{S}_j}^T)^T$ in each iteration for $1\leq k \leq d_n$ for each step.

We remark that Case 1 only utilises the confounding information from covariates $\mathbf{X}$ itself while Case 2 incorporates the exogenous conditional information but ignores the confounding effect from $\mathbf{X}$ itself. Case 3 is  the most flexible version incorporating all types of covariate information. We will implement Case 3 for the real data analysis in this paper.



\section{Numerical Studies}

\subsection{Simulation Studies}
In this section, we conduct simulations to examine the finite sample performances of the two proposed copula-based correlations: CC and CPC, as well as the two so-constructed screening procedures: CC-SIS and CPC-SIS.

\subsubsection{Inference Performance}
We consider two simulation examples with fixed dimension $p$ in this subsection, and illustrate the practical performance of estimated CC and CPC, respectively. We consider sample size $n=200$ and $400$ and set the number of repetitions to be $N=5000$ in Examples 1 and 2.

\textit{Example 1.} We generate the response from two models
\begin{itemize}
  \item (a1) $ Y = \exp(2X_1) + \exp((2+c_0)X_2)$, where $(X_1,X_2)$ is from a standard bivariate normal distribution with $corr(X_{1}, X_{2})=\rho$;
  \item (a2) $Y = 2X_{01} + (2+c_0)X_{02} + \varepsilon$,
\end{itemize}
where in model (a2), covariates $X_1$ and $X_2$ are both generated from a mixture distribution of a normal distribution with probability $0.9$ and a cauchy distribution with probability $0.1$, that is, $X_1 = 0.9X_{01}+0.1\epsilon_1$ and $X_2 = 0.9X_{02}+0.1\epsilon_2$ with $(X_{01}, X_{02})$ following a standard bivariate normal distribution with $corr(X_{01}, X_{02})=\rho$, $\epsilon_1\sim \frac{1}{5}Cauchy(0,1)$ and $\epsilon_2\sim \frac{1}{5}Cauchy(0,1)$ are independent and the model error $\varepsilon$ is generated from $N(0,1)$. Our interest of this example is to test $H_0: \varrho_{Y,X_1}(\tau, \iota)=\varrho_{Y,X_2}(\tau, \iota)$ at the significance level $\alpha_0=0.05$ under various values of $c_0$. We consider different $\rho$'s and set $c_0=0, 1,2,4$, where $c_0=0$ implies that $H_0$ holds true, whereas $H_0$ should be rejected with large probability for other values of $c_0$. We report the empirical size and power for each setup over 5000 runs in Table \ref{tab1}. Observing Table \ref{tab1}, we can see that our proposed CC testing procedure based on Theorem \ref{th1b} performs satisfactorily across different quantile levels $(\tau, \iota)$ since the values of empirical size are close to nominal level $\alpha_0=0.05$ and enlarging the sample size from $200$ to $400$ generally tends to improve the performance. Also we can see that when $c_0$ runs away from $0$, the empirical power increases to $1$, and when the correlation between $X_1$ and $X_2$ is low, the performance will be better.

\textit{Example 2.} In this example, we generate conditional variables $\mathbf{Z}=(Z_1, Z_2,Z_3, Z_4)^T$ from the 4-dimensional multivariate normal distribution $N(\mathbf{0}_4, \Sigma)$ with $\Sigma=(\rho^{|j-k|})_{1\leq j,k\leq 4}$. We generate $Y$ from the model
\begin{itemize}
  \item   $ Y = 2X_1 + (2+c_0)X_2 + \mathbf{Z}^T\mathbf{b} + \varepsilon,$
\end{itemize}
where $\mathbf{b}=(3,4,3,4)^T$, $X_1 = \mathbf{Z}^T\mathbf{b} + \epsilon_1 $ and $X_2 = \mathbf{Z}^T\mathbf{b} + \epsilon_2$, where $\epsilon_1\sim \frac{1}{3}t(3)$ and  $\epsilon_2\sim \frac{1}{3}t(3)$ are independent. The other setups are the same as Example 1. In this example, our interest is to test $H_0: \varrho_{Y,X_1|\mathbf{Z}}(\tau, \iota)=\varrho_{Y,X_2|\mathbf{Z}}(\tau, \iota)$ at the significance level $\alpha_0=0.05$. Here, we consider the performance of CPC with $(\tau,\iota)=(0.5,0.5)$ and the corresponding empirical size and power are reported in Table \ref{tab2}.  A similar conclusion can be drawn as in Example 1. These numerical results empirically demonstrate that our theoretical result in Theorem \ref{th3b} is valid.

\subsubsection{Screening Performance for CC-SIS and CPC-SIS}

Throughout this subsection, we adopt the following simulation setup: the sample size $n=200$, the covariate dimension $p_n=1000$, and the number of simulations $N=200$ for each simulation setup. Moreover, for the purpose of comparison, we use three criteria for evaluation: the first criterion is the minimum model size (MMS), namely, the smallest number of the selected covariates that contain all the active covariates, and its robust standard deviation (RSD); the second is the rank for each active covariates ($R_j$); and the third is the proportion of all the active covariates being selected ($\mathcal{P}$) with the screening threshold specified as $\lfloor n/\log n \rfloor$ over $N$ simulations. We report the median of MMS and $R_j$.

In Example 3, we compare our CC-SIS methods with a few existing methods: SIS (\cite{FL2008}), SIRS (\cite{ZLLZ2011}), DC-SIS (\cite{LZZ2012}), Kendall-SIS (\cite{LPZZ2012}), QC-SIS (\cite{LLT2015}) and CQC-SIS (\cite{MZ2016}). In Example 4, we compare our CPC-SIS procedure with the aforementioned marginal screening methods as well as the QPC-SIS (\cite{MLT2017}), where confounding effects arise from covariates $\mathbf{X}$. In this example, we employ the algorithm in Case 1 given in Section 4 in order to compare with the QPC-SIS by \cite{MLT2017}. In Example 5, we apply the algorithm in Case 2 and compare with the QPC-SIS by \cite{MLT2017}.

\textit{Example 3.} This example is used to assess the performance of the proposed CC-SIS. Let $\mathbf{X}^{\ast}=(X_1^{*},\ldots, X_{p_n}^{*})^T$ be a latent random vector having the  $p_n$-dimensional normal distribution $N(\mathbf{0}_{p_n}, \Sigma)$ with $\Sigma=(\rho^{|k-l|})_{1\leq k,l\leq p_n}$, where we set the correlation $\rho=0.4$ and $0.8$. We write $\ve\epsilon = (\epsilon_1,\ldots, \epsilon_{p_n})^T$ with each component $\epsilon_j$ being independent of other components and having the standard Cauchy distribution, i.e., $\epsilon_j\sim Cauchy(0,1)$. We generate covariates $\mathbf{X}$ from a mixture distribution: $\mathbf{X}=0.8\mathbf{X}^{*} + 0.2\ve\epsilon$, and simulate the response data from the following three models:
\begin{itemize}
  \item (b1) $ Y = 3X_1^{*} + 3X_2^{*} + 2X_3^{*} + 2X_{4}^{*} + 2X_{5}^{*} + \varepsilon$,
  \item (b2) $ Y = 5X_1^{*}I(X_1^{*}<0) + 5X_2^{*}I(X_2^{*}>0) + 5\sin(X_{10}^{*}) + \varepsilon$,
  \item (b3) $Y = \exp\{3\beta_1\sin(X_1^{*})+2\beta_2\exp(X_2^{*}) + 1.5\beta_3 I(X_3^{*}>0) + 2\log(|X_4^{*}|)\} + \varepsilon$,
\end{itemize}
where $\varepsilon$ is simulated from two scenarios: $\varepsilon\sim N(0,1)$ and $\varepsilon\sim Cauchy(0,1)$ and, in model (b3), we set $\beta_j=c_j(-1)^U(4\log n/\sqrt{n} + Z_0)$ for $j=1,2$ and $3$, where $U\sim Bernoulli(0.4)$, $Z_0\sim N(0,0.5^2)$ and $(c_1,c_2,c_3)=(1,0.5,1)$. The resulting screening results in terms of MMS and $\mathcal{P}$ are presented in Table \ref{tab3}.

Eyeballing Table \ref{tab3}, we can make some key observations. Under models (b1) and (b2), our CC-SIS outperforms SIS, SIRS, DC-SIS and QC-SIS. In this case, both response and covariates are heavy-tailed and thus traditional linear correlation screening methods all fail to work. Our methods are also comparable to the nonparametric Kendall's $\tau$ which achieves high accuracy but is slower due to the numerical integration in its implementation. (b3) is a difficult case and very hard to screen accurately. Under this case, our CC-SIS still performs much better than other methods.

\textit{Example 4.} This example is designed to evaluate the performance of the proposed CPC-SIS. We generate the response from the following two models
\begin{itemize}
  \item (d1) $ Y = \beta X_1 + \beta X_2 + \beta X_3 -3\beta \sqrt{\rho} X_{4} + \varepsilon$,
  \item (d2) $ Y = \beta X_1^{*} + \beta X_2^{*} + \beta X_3^{*} -3\beta \sqrt{\rho} X_{4}^{*} + \varepsilon$,
\end{itemize}
where $\beta=4$ for both models. In model (d1), we consider the observed covariates $\mathbf{X}\sim N(\mathbf{0}_{p_n}, \Sigma)$ with $\Sigma=(\sigma_{ij})_{1\leq i,j\leq p_n}$, implying that covariates $\mathbf{X}$ have no outliers and the response $Y$ is fully dependent on observed $\mathbf{X}$ up to a random noise. Under this setting, it is desired to expect that the QPC-SIS performs better than our CPC-SIS since $\mathbf{X}$ is normal. In model (d2), we generate covariates $\mathbf{X}$ from a mixture distribution $0.9\mathbf{X}^{*} + 0.1\ve\epsilon$, where $\mathbf{X}^{*} \sim N(\mathbf{0}_{p_n}, \Sigma)$ with $\Sigma=(\sigma_{ij})_{1\leq i,j\leq p_n}$, and each element of $\ve\epsilon$ is independent and distributed as $\frac{1}{5}Cauchy(0,1)$. In this model, the covariates are contaminated with outliers, while the heterogenicity
of response $Y$ stems merely from the random error $\varepsilon$. In addition, in model (d1), we let $\sigma_{ii}=1$ and $\sigma_{ij}=\rho, j\neq i$ except that $\sigma_{4j}=\sigma_{j4}=\sqrt{\rho}$. In model (d1), we set $\sigma_{ij}=0$ if $i>4$ or $j>4$ and the rest are the same as model (d1). For both models, at the population level, the covariate $X_4$ is marginally uncorrelated with $Y$. We consider two cases of $\rho=0.95$ and $0.5$ for simulation comparison.

Tables \ref{tab4a} and \ref{tab4b} report the screening results regarding the rank $R_j$ and MMS. From the tables, we can see that all the marginal screening approaches fail to pick out the covariate $X_4$ with very large values of the rank $R_4$. QPC-SIS and our CPC-SIS work much better than those marginal methods. Particularly, under model (d1), our CPC-SIS has a very competitive performance to QPC-SIS. Under model (d2), when the covariates are highly correlated ($\rho=0.95$), our proposal CPC-SIS$_{(0.5,0.5)}$ has the best performance.

\textit{Example 5.} In this example, we examine the case of $\mathbf{Z}=\mathbf{W}$. Since conditional variables selected for each $X_j$ are the same, so we can simply employ our proposed CPC-SIS and the QPC-SIS of \cite{MLT2017} for variable screening. We generate the response from the model
\begin{itemize}
  \item $ Y = 2 X_1 + 2 X_2 - 4 X_3 +3 X_{4} + \varepsilon$,
\end{itemize}
where $X_j=\mathbf{W}^T\mathbf{b}+U_j$, $\mathbf{W}$ has the same distribution as $\mathbf{Z}$ in Example 2, $\mathbf{b}=(2,4/3,2, 4/3)^T$ and $U_j\sim \frac{1}{3}Cauchy(0,1)$ for $j=1,\cdots, p_n$. The model error $\varepsilon$ is the same as that in Example 4. The simulation results are given in Table \ref{tab5}. As expected, we can see that all the marginal screening procedures fail to work since they are unable to identify the covariate $X_3$. Our proposed CPC-SIS still performs better than QPC-SIS in terms of MMS.

\subsection{Real Data Applications}

\subsubsection{Rats Data}
\begin{figure}
  \centering
  \includegraphics[width=6.5cm, height=6cm]{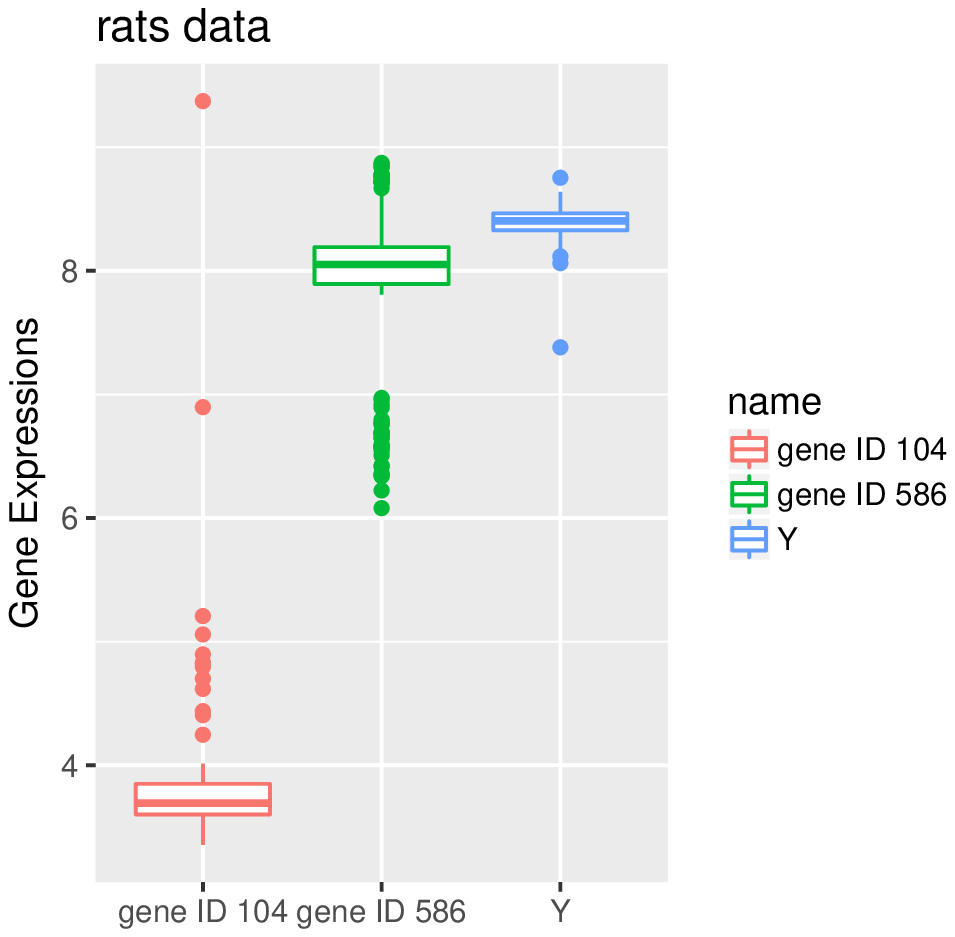}
  \includegraphics[width=6.5cm, height=6cm]{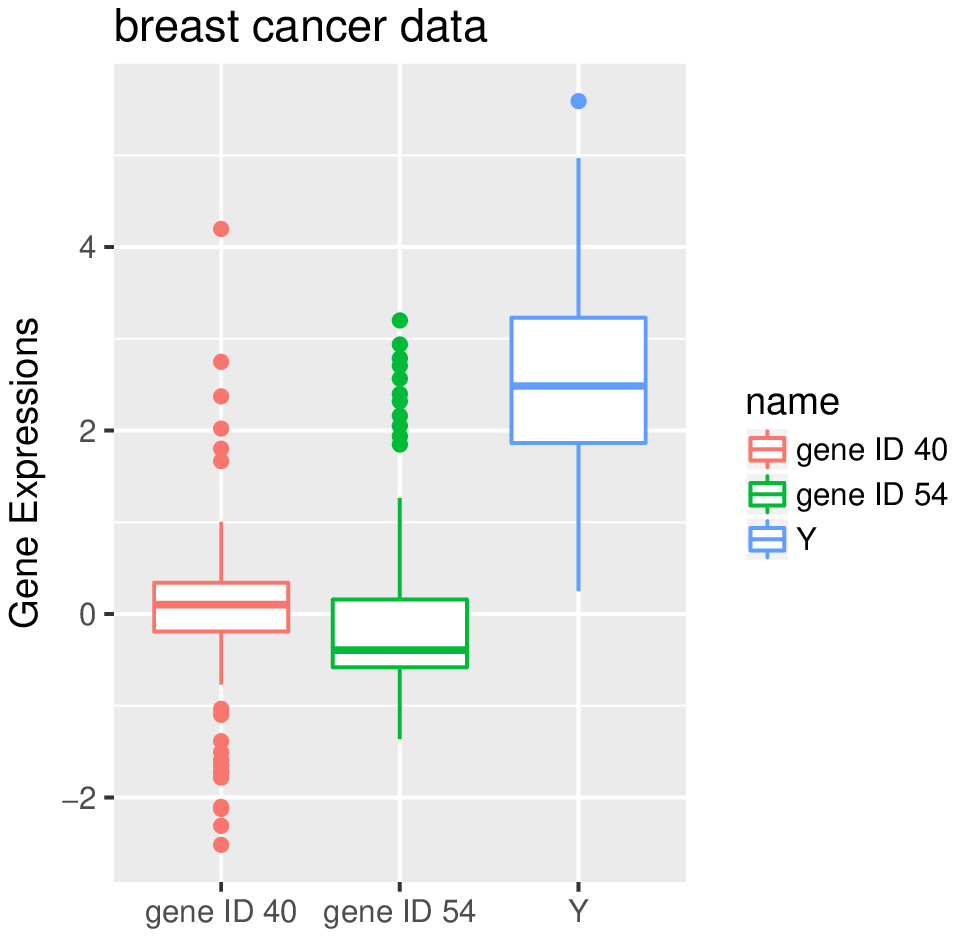}
  \caption{ Box-plots for the response and two randomly selected genes for the two datasets. The left panel is for the rats data and the right panel is for the breast cancer data.  }\label{fig1}
\end{figure}

We illustrate the CC-SIS and CPC-SIS with the gene expression data on 120 male rats of 12 weeks old, including expression measurements of 31,099 gene probes. It has been analysed in \cite{Set2006} for investigation of the gene regulation in the mammalian and is available at \url{ftp://ftp.ncbi.nlm.nih.gov/geo/series/GSE5nnn/GSE5680/matrix}. We follow \cite{MLT2017} and consider the expression of gene TRIM32 (probe 1389163\_at) as the response variable $Y$ since it was identified to cause Bardet-Biedl syndrome, closely associated with the human hereditary disease of the retina \cite{Cet2006}. The other gene probes are treated as the covariates ${\bf X}$. We first apply \cite{IH1993}'s (IH) approach to check outliers. IH constructs a Z-score $Z_i=0.6745(x_i-\tilde{x})/\mathrm{MAD}$, where $\mathrm{MAD}$ denotes the median absolute deviation and $\tilde{x}$ stands for the median, and recommends that any $i$ such that $Z_i >3.5$ be labeled as outliers. IH method is quite popular in real applications such as engineering.  Consequently we find that there are over $60\%$ gene probes having one or more outliers. Figure \ref{fig1} displays the box-plots for two selected genes and the the response. If one only employs the conventional screening method ignoring the outliers, it would lead to inappropriate results. The copula-based methods may thus be more robust in this situation. In this data analysis, we have the sample $\{(Y_i,\mathbf{X}_i\in \mathbb{R}^{p_n}), 1\leq i\leq n\}$ with $n=120$ and $p_n=31098$.

 We report the overlaps of the top $\lfloor n/\log n \rfloor=25$ selected genes by various methods in Table \ref{tab:d2}. We can see that different methods select quite different genes and such low level of agreement should not be overlooked in practice. Robust and joint screening methods like what we propose in this paper lead to entirely different set of genes which are otherwise screened out by the conventional non-robust and marginal screening approaches. We notice that the CPC-SIS and QPC-SIS have a couple of overlaps, partly because both are conditional screening procedures and able to adjust the confounder effects.

 Table \ref{tab:d3} gives a summary of top 10 gene probes by different methods along with the p-value resulted from a marginal Wald-test. We then use these 10 genes as regressors and build a joint statistical models to predict $Y$. Linear regression and quantile regression are both considered for this purpose and we \textcolor[rgb]{1.00,0.00,0.00}{display} the mean of their prediction errors (PE1 and PE2) over 500 random partitions, where the partition ratio of training sample to test sample is $4:1$ for each partition. The PE is computed as the average of $\{(Y_i-\widehat{Y}_i)^2, i\in \text{testing set}\}$ and $\widehat{Y}_i$ is the predicted value at the $i$th test data point using the model constructed by the training sample with the 10 genes in Table \ref{tab:d3}. We can see that our proposed copula-based partial correlation screening performs the best with the smallest prediction error. Such a nice prediction result may be attributed to the fact that CPC selects appropriate markers for joint modelling after addressing the distribution heterogeneity and the conditional effects. The heterogeneity problem typically inflates the variance while a purely marginal screener could introduce bias. The prediction error, consists of the variance and the bias components, is thus much smaller after employing the CPC screening method.


\subsubsection{Breast Cancer Data}
The second data we use to illustrate our proposal is breast cancer data. Breast cancer has become the second most common cancer in the world and the most leading cause in women. There were nearly 1.7 million new cases diagnosed in 2012, according to the worldwide statistics\footnote{\url{http://www.wcrf.org/int/cancer-facts-figures/worldwide-data}}. Meanwhile, approximately 252,710 new cases of invasive breast cancer and 40,610 breast cancer deaths are expected to occur among US women in 2017 as reported in \cite{Det2017}.   Although major progresses in breast cancer treatment were made, there is also limited ability to predict the metastatic behavior of tumor. \cite{Vet2002} was the first to study the breast cancer study involving 97 lymph node-negative breast cancer patients 55 years old or younger, of which 46 developed distant metastases within 5 years (metastatic outcome coded as 1) and 51 remained metastases free for at least 5 years (metastatic outcome coded as 0). This expression data set with clinical variables has been well analysed in many papers for classification (\cite{BPD2008}, \cite{YLM2012}, among others).

In this study, after removing the genes with missing values, there are expression levels of 24,188 gene probes entering into the next analysis. In addition to gene expression measurements, the data for several clinical factors are available as well. Our interest is to identify which gene probes affect the tumor size given other clinical factors $(\mathbf{W})$ including age, histological grade, angioinvasion, lymphocytic infiltration, estrogen receptor and progesterone receiptor status. Therefore, we have the data $\{(Y_i,\mathbf{X}_i\in \mathbb{R}^{p_n}, \mathbf{W}_i\in \mathbb{R}^{r}), 1\leq i\leq n\}$ with $n=97$, $p_n=24,188$ and $r=6$ for further analysis.

Using the IH method on outlier detection, we find that 18,098 gene probes have at least 1 and at most 29 outliers, suggesting that approximately three quarters of overall gene probes contain extremely large values. The right panel of Figure \ref{fig1} displays the empirical distribution of the response and two typical covariates. Thus,  it is more suitable to apply robust joint screening approach such as the proposed CPC-SIS. We consider the three cases discussed in Section 4 and denote the methods as CPC-SISa$_1$, CPC-SISa$_2$ and CPC-SISa$_3$, respectively). Table \ref{tab:d4} gives the overlaps of the selected genes by various methods. A similar conclusion to that in the rats data analysis can be made. Furthermore, Table \ref{tab:d5} presents a summary of top 10 gene probes selected by various methods. The results on PE1 and PE2 in Table \ref{tab:d5} empirically verifies that our proposed CPC-SIS in Case 3 has the most satisfactory performance in out-of-sample prediction.

\section{Conclusion and Discussion}
We propose a copula-based correlation and partial correlation to facilitate robust marginal and joint screening for ultrahigh dimensional data sets. Large sample properties for the estimated correlation and sure screening properties for CC and CPC screeners were provided. Empirical studies including simulations and two data applications show that our proposed CC-SIS and CPC-SIS outperform the existing variable screening approaches, when outliers are present in both covariates and response. Therefore, our current proposals are more applicable to the ultrahigh dimensional heterogeneous data. We provide a guideline to carry out variable screening as follows. If the response and predictors are all normal without heteroscedastic variance and the predictors have low correlation, any marginal screening methods (SIS, SIRS, DC-SIS) can be applied. If the response contains outliers or follows a heavy tail distribution and the covariates are normal, robust screening methods (Kendall SIS, QC-SIS, CQC-SIS, CC-SIS) can be employed. If the covariates are highly correlated and conditional variables are available, conditional screening procedures (CSIS, QPC-SIS) can be used. If the data is heteroscedastic for both the response and covariates and the covariates may be highly correlated, then only CPC-SIS can be recommended.

The copula formulation may suggest many possible extensions of our methodology. Firstly, we may consider censored survival time outcome in this framework. See \cite{YL2018,HL2018,HMQ2018} for recent reviews on feature selection and screening for survival analysis. The estimation of copula-based correlation and partial correlation needs to incorporate the random censoring for such data and we need to invoke more complicated empirical process theories to argue the weak convergence results. Secondly, we may even allow the predictors to be censored. See \cite{CF2008} and \cite{CL2015} for some earlier discussion. Thirdly, we may consider more pairs of $(\tau, \iota)$ over a candidate set or an interval to incorporate more information on quantile of response and covariates. The relevant theoretical results in this paper can be further generalised. Such extensions to ultrahigh-dimensional data is non-trivial and requires a further development in the future work.


\noindent\textbf{Acknowledgement} We thank the Associate Editor and two anonymous referees for constructive comments on our manuscript. This work was partially supported by National Natural Science Foundation of China (Grant No. 11801202), Academic Research Fundings from Ministry of Eduction (MOE) in Singapore: R-155-000-197-112, R-155-000-195-114 and R-155-000-205-114.

\smallskip
\noindent\textbf{Supplementary Materials} The supplementary materials consists of two parts: Appendix A and Appendix B. The proofs of all theoretical results stated in the manuscript are given in the Appendix A. Some additional simulations are considered in the Appendix B, which includes the empirical performance of CC and CPC estimators and some discussion on the choice of quantile levels $(\tau,\iota)$ involved in the proposed approaches.


\newpage

\begin{table}[htbp]
 {  \footnotesize
\tabcolsep 0pt \vspace*{-12pt}
\def\temptablewidth{1.0\textwidth}
\centering
  \caption{\small Empirical size and power for testing $H_0: \varrho_{Y,X_1}(\tau, \iota)=\varrho_{Y,X_2}(\tau, \iota)$ using CC across $(\tau,\iota)=(0.25, 0.25)$, $(0.5, 0.5)$ and $(0.75, 0.75)$ for Example 1.  }  \label{tab1}%
 \begin{tabular*}{\temptablewidth}{@{\extracolsep{\fill}} llc cccccccc}
\hline\hline
      &       &       & \multicolumn{4}{c|}{Model (a1)} & \multicolumn{4}{c}{Model (a2)} \\
\cline{4-11}$\rho$ & \multicolumn{1}{l}{Method$(\tau,\iota)$} & $n$   & $c_0= 0$ & 1     & 2     & 4     & $c_0= 0$ & 1     & 2     & 4 \\
\hline
$0$ & CC$_{(0.25, 0.25)}$ & 200   & 0.056 & 0.07  & 0.07  & 0.08  & 0.051 & 0.48  & 0.90  & 1.00 \\
      &       & 400   & 0.051 & 0.07  & 0.09  & 0.08  & 0.055 & 0.78  & 1.00  & 1.00 \\
      & CC$_{(0.5, 0.5)}$ & 200   & 0.059 & 0.33  & 0.69  & 0.94  & 0.057 & 0.56  & 0.95  & 1.00 \\
      &       & 400   & 0.054 & 0.57  & 0.93  & 1.00  & 0.048 & 0.85  & 1.00  & 1.00 \\
      & CC$_{(0.75, 0.75)}$ & 200   & 0.057 & 0.77  & 0.99  & 1.00  & 0.053 & 0.49  & 0.90  & 1.00 \\
      &       & 400   & 0.055 & 0.96  & 1.00  & 1.00  & 0.050 & 0.78  & 1.00  & 1.00 \\   \hline
$0.5$ & CC$_{(0.25, 0.25)}$ & 200   & 0.055 & 0.06  & 0.08  & 0.14  & 0.048 & 0.27  & 0.62  & 0.94 \\
      &       & 400   & 0.056 & 0.06  & 0.08  & 0.18  & 0.052 & 0.45  & 0.90  & 1.00 \\
      & CC$_{(0.5, 0.5)}$ & 200   & 0.057 & 0.26  & 0.56  & 0.85  & 0.050 & 0.29  & 0.68  & 0.96 \\
      &       & 400   & 0.053 & 0.44  & 0.83  & 0.98  & 0.049 & 0.51  & 0.93  & 1.00 \\
      & CC$_{(0.75, 0.75)}$ & 200   & 0.054 & 0.71  & 0.99  & 1.00  & 0.052 & 0.26  & 0.63  & 0.93 \\
      &       & 400   & 0.052 & 0.95  & 1.00  & 1.00  & 0.043 & 0.47  & 0.89  & 1.00 \\  \hline
$0.9$ & CC$_{(0.25, 0.25)}$ & 200   & 0.055 & 0.08  & 0.17  & 0.47  & 0.048 & 0.10  & 0.19  & 0.43 \\
      &       & 400   & 0.052 & 0.09  & 0.29  & 0.75  & 0.052 & 0.14  & 0.35  & 0.73 \\
      & CC$_{(0.5, 0.5)}$ & 200   & 0.061 & 0.17  & 0.36  & 0.63  & 0.048 & 0.09  & 0.20  & 0.45 \\
      &       & 400   & 0.050 & 0.27  & 0.59  & 0.87  & 0.048 & 0.15  & 0.36  & 0.75 \\
      & CC$_{(0.75, 0.75)}$ & 200   & 0.045 & 0.52  & 0.91  & 0.99  & 0.053 & 0.09  & 0.20  & 0.43 \\
      &       & 400   & 0.047 & 0.85  & 1.00  & 1.00  & 0.049 & 0.14  & 0.36  & 0.72 \\
\hline \hline
    \end{tabular*}%
}
\end{table}

\begin{table}[htbp]
 { \footnotesize
\tabcolsep 0pt \vspace*{-12pt}
\def\temptablewidth{1.0\textwidth}
\centering
  \caption{\small Empirical size and power for testing $H_0: \varrho_{Y,X_1|\mathbf{Z}}(\tau, \iota)=\varrho_{Y,X_2|\mathbf{Z}}(\tau, \iota)$ using CPC with $(\tau,\iota)=(0.5, 0.5)$ in Example 2.  }  \label{tab2}%
\begin{tabular*}{\temptablewidth}{@{\extracolsep{\fill}} lc cccccccc}
\hline\hline
      &             & \multicolumn{4}{c}{$\varepsilon\sim N(0, 1)$} & \multicolumn{4}{c}{$\varepsilon\sim \frac{1}{3}Cauchy(0, 1)$} \\ \cline{3-6}\cline{7-10}
$\rho$  & $n$   & $c_0= 0$ & 1     & 2     & 4     & $c_0= 0$ & 1     & 2     & 4 \\
\hline
$0$  & $200$ & 0.074 & 0.40 & 0.84 & 1.00 & 0.084 & 0.45 & 0.87 & 1.00  \\
     & $400$ & 0.067 & 0.65 & 0.98 & 1.00 & 0.073 & 0.72 & 0.99 & 1.00 \\
$0.5$ & $200$ & 0.073 & 0.40 & 0.83 & 1.00 & 0.081 & 0.46 & 0.86 & 1.00 \\
     & $400$ & 0.066 & 0.65 & 0.98 & 1.00 & 0.068 & 0.72 & 0.99 & 1.00 \\
$0.9$ & $200$ & 0.077 & 0.40 & 0.84 & 1.00 & 0.080 & 0.46 & 0.86 & 1.00 \\
     & $400$ & 0.069 & 0.65 & 0.98 & 1.00 & 0.070 & 0.72 & 0.99 & 1.00 \\
\hline \hline
\end{tabular*}%
}
\end{table}

\begin{table}[htbp]
 {\footnotesize 
\tabcolsep 0pt \vspace*{-12pt}
\def\temptablewidth{1.0\textwidth}
\centering
  \caption{\small Simulation results for Example 3, where MMS stands for the median of the minimum model size and its robust standard deviations (RSD) are given in parenthesis, $\mathcal{P}$ is the proportion of screened sets that cover all active predictors with screening parameter $d_n= \lfloor n/\log n \rfloor$. }  \label{tab3}%
\begin{tabular*}{\temptablewidth}{@{\extracolsep{\fill}} l cc cc cc cc}
    \hline \hline
	&	\multicolumn{4}{c}{$\rho=0.4$}			&	\multicolumn{4}{c}{$\rho=0.8$}									\\ \cline{2-5}\cline{6-9}
	&	\multicolumn{2}{c}{$\varepsilon\sim N(0,1)$}		&	\multicolumn{2}{c}{$\varepsilon\sim Cauchy(0,1)$}	&	\multicolumn{2}{c}{$\varepsilon\sim N(0,1)$}	&	\multicolumn{2}{c}{$\varepsilon\sim Cauchy(0,1)$}		\\ \cline{2-3}\cline{4-5}\cline{6-7}\cline{8-9}
Method$(\tau, \iota)$  &	MMS(RSD)	&	$\mathcal{P}$ &	MMS(RSD)	&	$\mathcal{P}$  &	MMS(RSD)	&	$\mathcal{P}$  &	MMS(RSD)	&	$\mathcal{P}$   \\ \hline
          \multicolumn{9}{c}{Model (b1)} \\
SIS        & 655        (396)      & 0.03       & 790        (227)      & 0.00       & 522        (431)      & 0.12       & 608        (389)      & 0.03 \\
SIRS       & 615        (176)      & 0.00       & 669        (158)      & 0.00       & 513        (149)      & 0.00       & 547        (138)      & 0.00 \\
DC-SIS     & 634        (441)      & 0.09       & 725        (350)      & 0.04       & 460        (527)      & 0.20       & 615        (489)      & 0.13 \\
Kendall-SIS    & 5          (0)        & 0.99       & 5          (1)        & 0.96       & 5          (0)        & 1.00       & 5          (0)        & 1.00 \\
CC-SIS$_{(0.25,0.25)}$ & 9          (18)       & 0.80       & 20         (41)       & 0.70       & 5          (0)        & 1.00       & 5          (0)        & 1.00 \\
CC-SIS$_{(0.5,0.5)}$ & 6          (7)        & 0.93       & 9          (14)       & 0.84       & 5          (0)        & 1.00       & 5          (0)        & 1.00 \\
CC-SIS$_{(0.75,0.75)}$ & 9          (15)       & 0.85       & 14         (40)       & 0.70       & 5          (0)        & 1.00       & 5          (0)        & 1.00 \\
QC-SIS$_{(0.25)}$ & 710        (315)      & 0.01       & 782        (271)      & 0.00       & 649        (444)      & 0.06       & 678        (406)      & 0.03 \\
QC-SIS$_{(0.5)}$ & 716        (334)      & 0.01       & 790        (250)      & 0.02       & 665        (438)      & 0.07       & 693        (398)      & 0.06 \\
QC-SIS$_{(0.75)}$ & 711        (297)      & 0.01       & 751        (314)      & 0.01       & 599        (379)      & 0.04       & 670        (408)      & 0.02 \\
CQC-SIS    & 573        (156)      & 0.00       & 614        (149)      & 0.00       & 442        (118)      & 0.00       & 474        (104)      & 0.00 \\
\hline
           \multicolumn{9}{c}{Model (b2)}  \\
SIS        & 482        (408)      & 0.13       & 638        (376)      & 0.01       & 351        (414)      & 0.15       & 686        (331)      & 0.03 \\
SIRS       & 551        (161)      & 0.00       & 589        (195)      & 0.00       & 523        (173)      & 0.00       & 543        (158)      & 0.00 \\
DC-SIS     & 284        (585)      & 0.30       & 607        (493)      & 0.12       & 230        (569)      & 0.35       & 506        (514)      & 0.17 \\
Kendall-SIS    & 3          (0)        & 1.00       & 3          (0)        & 1.00       & 4          (1)        & 1.00       & 4          (1)        & 1.00 \\
CC-SIS$_{(0.25,0.25)}$ & 6          (20)       & 0.78       & 13         (45)       & 0.68       & 5          (2)        & 0.99       & 7          (4)        & 0.94 \\
CC-SIS$_{(0.5,0.5)}$ & 3          (1)        & 0.98       & 4          (3)        & 0.95       & 4          (2)        & 1.00       & 5          (1)        & 1.00 \\
CC-SIS$_{(0.75,0.75)}$ & 7          (15)       & 0.85       & 14         (43)       & 0.66       & 6          (3)        & 0.99       & 7          (3)        & 0.96 \\
QC-SIS$_{(0.25)}$ & 615        (356)      & 0.06       & 588        (391)      & 0.05       & 411        (418)      & 0.07       & 598        (411)      & 0.05 \\
QC-SIS$_{(0.5)}$ & 584        (460)      & 0.07       & 640        (373)      & 0.05       & 539        (474)      & 0.13       & 592        (425)      & 0.10 \\
QC-SIS$_{(0.75)}$ & 588        (414)      & 0.03       & 579        (353)      & 0.02       & 486        (374)      & 0.09       & 604        (401)      & 0.06 \\
CQC-SIS    & 512        (153)      & 0.00       & 562        (155)      & 0.00       & 441        (122)      & 0.00       & 481        (150)      & 0.00 \\
\hline
           \multicolumn{9}{c}{Model (b3)}  \\
SIS        & 737        (296)      & 0.00       & 768        (204)      & 0.00       & 524        (399)      & 0.01       & 619        (368)      & 0.01 \\
SIRS       & 767        (215)      & 0.00       & 745        (182)      & 0.00       & 664        (240)      & 0.00       & 703        (203)      & 0.00 \\
DC-SIS     & 790        (245)      & 0.00       & 789        (213)      & 0.01       & 722        (276)      & 0.02       & 719        (309)      & 0.01 \\
Kendall-SIS    & 504        (467)      & 0.12       & 475        (442)      & 0.10       & 5          (192)      & 0.65       & 16         (377)      & 0.60 \\
CC-SIS$_{(0.25,0.25)}$ & 475        (234)      & 0.02       & 679        (317)      & 0.02       & 252        (365)      & 0.21       & 279        (455)      & 0.14 \\
CC-SIS$_{(0.5,0.5)}$ & 483        (357)      & 0.09       & 478        (337)      & 0.07       & 18         (341)      & 0.56       & 53         (357)      & 0.47 \\
CC-SIS$_{(0.75,0.75)}$ & 426        (419)      & 0.17       & 281        (460)      & 0.19       & 12         (309)      & 0.59       & 113        (349)      & 0.47 \\
QC-SIS$_{(0.25)}$ & 791        (199)      & 0.00       & 795        (191)      & 0.00       & 794        (222)      & 0.01       & 781        (225)      & 0.01 \\
QC-SIS$_{(0.5)}$ & 785        (210)      & 0.00       & 810        (209)      & 0.01       & 750        (266)      & 0.02       & 796        (219)      & 0.04 \\
QC-SIS$_{(0.75)}$ & 762        (257)      & 0.00       & 793        (240)      & 0.00       & 698        (408)      & 0.01       & 714        (311)      & 0.02 \\
CQC-SIS    & 743        (216)      & 0.00       & 749        (177)      & 0.00       & 623        (210)      & 0.00       & 628        (260)      & 0.00 \\
\hline \hline
\end{tabular*}%
}
\end{table}

\begin{table}
{ \footnotesize  
\tabcolsep 0pt \vspace*{-12pt}
\def\temptablewidth{1.0\textwidth}
\centering
\caption{\small Simulation results for Example 4 (d1), where $R_j$ indicates the median of the rank of the relevant predictors and MMS stands for the median of the minimum model size and its robust standard deviations (RSD) are given in parenthesis. }  \label{tab4a}%
\begin{tabular*}{\temptablewidth}{@{\extracolsep{\fill}} ll ccccc ccccc}
\hline\hline
      &       & \multicolumn{5}{c}{$\varepsilon\sim N(0,1)$}                    & \multicolumn{5}{c}{$\varepsilon\sim \frac{1}{3}Cauchy(0,1)$}   \\ \cline{3-7}\cline{8-12}
 $\rho$ &  Method$(\tau, \iota)$     & $R_1$    & $R_2$    & $R_3$    & $R_4$    & MMS (RSD) & $R_1$    & $R_2$    & $R_3$    & $R_4$    & MMS (RSD) \\
\hline
\multicolumn{12}{c}{Model (d1)} \\
0.5   & SIS   & 2     & 2     & 2     & 450   & 450  (235) & 3     & 3     & 3     & 411   & 476  (270) \\
      & SIRS  & 3     & 3     & 3     & 469   & 501  (243) & 3     & 3     & 4     & 432   & 466  (258) \\
      & DC-SIS & 2     & 2     & 2     & 497   & 497  (280) & 2     & 2     & 2     & 464   & 464  (220) \\
      & Kendall-SIS & 2     & 2     & 2     & 454   & 454  (228) & 2     & 2     & 2     & 467   & 468  (280) \\
      & CC-SIS$_{(0.25,0.25)}$ & 2     & 2     & 2     & 367   & 372  (331) & 2     & 3     & 2     & 382   & 390  (311) \\
      & CC-SIS$_{(0.5,0.5)}$ & 2     & 2     & 2     & 415   & 425  (272) & 2     & 2     & 2     & 396   & 408  (275) \\
      & CC-SIS$_{(0.75,0.75)}$ & 2     & 2     & 2     & 371   & 387  (324) & 2     & 3     & 2     & 378   & 388  (311) \\
      & QC-SIS$_{(0.25)}$ & 2     & 2     & 2     & 478   & 493  (270) & 2     & 2     & 2     & 473   & 477  (245) \\
      & QC-SIS$_{(0.5)}$ & 2     & 2     & 2     & 437   & 437  (236) & 2     & 2     & 2     & 461   & 463  (331) \\
      & QC-SIS$_{(0.75)}$ & 2     & 2     & 2     & 437   & 438  (227) & 2     & 2     & 2     & 408   & 421  (247) \\
      & QPC-SIS$_{(0.25)}$ & 2     & 2     & 2     & 4     & 4    (0)   & 2     & 2     & 2     & 4     & 4    (0) \\
      & QPC-SIS$_{(0.5)}$ & 2     & 2     & 2     & 4     & 4    (0)   & 2     & 2     & 2     & 4     & 4    (0) \\
      & QPC-SIS$_{(0.75)}$ & 2     & 2     & 2     & 4     & 4    (0)   & 2     & 2     & 2     & 4     & 4    (0) \\
      & CPC-SIS$_{(0.25,0.25)}$ & 2     & 2     & 2     & 4     & 4    (0)   & 2     & 2     & 2     & 4     & 4    (0) \\
      & CPC-SIS$_{(0.5,0.5)}$ & 2     & 2     & 2     & 4     & 4    (0)   & 2     & 2     & 2     & 4     & 4    (0) \\
      & CPC-SIS$_{(0.75,0.75)}$ & 2     & 2     & 2     & 4     & 4    (0)   & 2     & 2     & 2     & 4     & 4    (0) \\   \hline
0.95  & SIS   & 3     & 3     & 3     & 470   & 498  (306) & 214   & 269   & 262   & 466   & 735  (353) \\
      & SIRS  & 28    & 52    & 45    & 474   & 541  (331) & 24    & 26    & 43    & 473   & 584  (376) \\
      & DC-SIS & 3     & 3     & 3     & 506   & 556  (220) & 10    & 12    & 15    & 495   & 610  (240) \\
      & Kendall-SIS & 3     & 3     & 3     & 482   & 529  (296) & 3     & 3     & 3     & 486   & 557  (325) \\
      & CC-SIS$_{(0.25,0.25)}$ & 91    & 95    & 104   & 303   & 481  (276) & 88    & 99    & 88    & 318   & 485  (292) \\
      & CC-SIS$_{(0.5,0.5)}$ & 56    & 49    & 67    & 365   & 499  (383) & 53    & 68    & 47    & 354   & 610  (377) \\
      & CC-SIS$_{(0.75,0.75)}$ & 114   & 92    & 94    & 297   & 470  (301) & 125   & 63    & 138   & 318   & 492  (289) \\
      & QC-SIS$_{(0.25)}$ & 3     & 3     & 3     & 470   & 512  (387) & 3     & 3     & 3     & 467   & 522  (375) \\
      & QC-SIS$_{(0.5)}$ & 3     & 3     & 3     & 465   & 499  (346) & 3     & 3     & 3     & 467   & 506  (382) \\
      & QC-SIS$_{(0.75)}$ & 3     & 3     & 3     & 468   & 520  (379) & 3     & 3     & 3     & 469   & 508  (329) \\
      & QPC-SIS$_{(0.25)}$ & 2     & 2     & 2     & 4     & 4    (0)   & 2     & 2     & 2     & 4     & 4    (0) \\
      & QPC-SIS$_{(0.5)}$ & 2     & 2     & 2     & 4     & 4    (0)   & 2     & 2     & 2     & 4     & 4    (0) \\
      & QPC-SIS$_{(0.75)}$ & 2     & 2     & 2     & 4     & 4    (0)   & 2     & 2     & 2     & 4     & 4    (0) \\
      & CPC-SIS$_{(0.25,0.25)}$ & 2     & 2     & 2     & 4     & 4    (1)   & 2     & 2     & 2     & 5     & 5    (1) \\
      & CPC-SIS$_{(0.5,0.5)}$ & 2     & 2     & 2     & 4     & 4    (0)   & 2     & 2     & 2     & 4     & 4    (0) \\
      & CPC-SIS$_{(0.75,0.75)}$ & 2     & 2     & 2     & 4     & 4    (1)   & 2     & 2     & 2     & 4     & 4    (1) \\
\hline\hline
\end{tabular*}%
}
\end{table}%

\begin{table}
{\footnotesize  
\tabcolsep 0pt \vspace*{-12pt}
\def\temptablewidth{1.0\textwidth}
\centering
\caption{\small Simulation results for Example 4 (d2), where $R_j$ indicates the median of the rank of the relevant predictors and MMS stands for the median of the minimum model size and its robust standard deviations (RSD) are given in parenthesis. }  \label{tab4b}%
\begin{tabular*}{\temptablewidth}{@{\extracolsep{\fill}} ll ccccc ccccc}
\hline\hline
      &       & \multicolumn{5}{c}{$\varepsilon\sim N(0,1)$}                    & \multicolumn{5}{c}{$\varepsilon\sim \frac{1}{3}Cauchy(0,1)$}   \\ \cline{3-7}\cline{8-12}
 $\rho$ &  Method$(\tau, \iota)$     & $R_1$    & $R_2$    & $R_3$    & $R_4$    & MMS (RSD) & $R_1$    & $R_2$    & $R_3$    & $R_4$    & MMS (RSD) \\
\hline
\multicolumn{12}{c}{Model (d2)} \\
0.5   & SIS   & 2     & 2     & 2     & 530   & 566  (359) & 10    & 12    & 14    & 535   & 667  (351) \\
      & SIRS  & 99    & 93    & 100   & 453   & 485  (304) & 107   & 110   & 112   & 478   & 507  (356) \\
      & DC-SIS & 2     & 2     & 2     & 558   & 629  (393) & 2     & 2     & 2     & 554   & 638  (322) \\
      & Kendall-SIS & 2     & 2     & 2     & 514   & 514  (362) & 2     & 2     & 2     & 578   & 578  (381) \\
      & CC-SIS$_{(0.25,0.25)}$ & 3     & 3     & 3     & 437   & 442  (362) & 4     & 4     & 4     & 453   & 456  (333) \\
      & CC-SIS$_{(0.5,0.5)}$ & 3     & 3     & 3     & 339   & 348  (358) & 2     & 3     & 3     & 459   & 463  (343) \\
      & CC-SIS$_{(0.75,0.75)}$ & 3     & 3     & 3     & 440   & 446  (315) & 3     & 4     & 5     & 449   & 451  (333) \\
      & QC-SIS$_{(0.25)}$ & 3     & 3     & 3     & 499   & 594  (352) & 3     & 3     & 2     & 531   & 618  (363) \\
      & QC-SIS$_{(0.5)}$ & 2     & 2     & 2     & 514   & 583  (377) & 2     & 3     & 2     & 501   & 533  (322) \\
      & QC-SIS$_{(0.75)}$ & 3     & 3     & 3     & 532   & 588  (345) & 3     & 3     & 3     & 460   & 537  (337) \\
      & QPC-SIS$_{(0.25)}$ & 3     & 3     & 3     & 2     & 4    (1)   & 3     & 3     & 3     & 2     & 4    (1) \\
      & QPC-SIS$_{(0.5)}$ & 2     & 3     & 3     & 2     & 4    (0)   & 3     & 3     & 3     & 2     & 4    (0) \\
      & QPC-SIS$_{(0.75)}$ & 3     & 3     & 3     & 2     & 4    (1)   & 3     & 3     & 3     & 2     & 4    (3) \\
      & CPC-SIS$_{(0.25,0.25)}$ & 2     & 2     & 2     & 4     & 4    (2)   & 2     & 2     & 2     & 4     & 5    (3) \\
      & CPC-SIS$_{(0.5,0.5)}$ & 3     & 3     & 3     & 1     & 4    (0)   & 3     & 3     & 3     & 2     & 4    (0) \\
      & CPC-SIS$_{(0.75,0.75)}$ & 2     & 2     & 2     & 4     & 4    (1)   & 2     & 2     & 3     & 4     & 5    (3) \\
\hline
0.95  & SIS   & 182   & 169   & 200   & 557   & 712  (285) & 518   & 434   & 435   & 496   & 770  (282)\\
      & SIRS  & 355   & 388   & 375   & 507   & 729  (290) & 450   & 401   & 403   & 507   & 720  (305) \\
      & DC-SIS & 238   & 191   & 236   & 540   & 762  (335) & 347   & 329   & 365   & 507   & 825  (262) \\
      & Kendall-SIS & 192   & 154   & 181   & 534   & 675  (274) & 224   & 222   & 196   & 482   & 734  (309) \\
      & CC-SIS$_{(0.25,0.25)}$ & 266   & 281   & 268   & 443   & 695  (197) & 254   & 260   & 255   & 438   & 688  (205) \\
      & CC-SIS$_{(0.5,0.5)}$ & 292   & 224   & 303   & 473   & 672  (305) & 318   & 305   & 213   & 462   & 666  (300) \\
      & CC-SIS$_{(0.75,0.75)}$ & 262   & 265   & 265   & 441   & 694  (199) & 268   & 423   & 284   & 431   & 693  (195) \\
      & QC-SIS$_{(0.25)}$ & 313   & 295   & 338   & 595   & 791  (259) & 393   & 394   & 357   & 480   & 755  (280) \\
      & QC-SIS$_{(0.5)}$ & 274   & 304   & 308   & 529   & 795  (283) & 291   & 282   & 290   & 480   & 753  (258) \\
      & QC-SIS$_{(0.75)}$ & 259   & 341   & 279   & 529   & 746  (276) & 362   & 374   & 376   & 459   & 806  (268) \\
      & QPC-SIS$_{(0.25)}$ & 112   & 150   & 158   & 4     & 530  (465) & 163   & 186   & 162   & 5     & 583  (475) \\
      & QPC-SIS$_{(0.5)}$ & 106   & 147   & 103   & 3     & 518  (524) & 148   & 102   & 96    & 3     & 514  (535) \\
      & QPC-SIS$_{(0.75)}$ & 148   & 164   & 214   & 3     & 574  (437) & 146   & 196   & 126   & 6     & 591  (483) \\
      & CPC-SIS$_{(0.25,0.25)}$ & 14    & 11    & 11    & 3     & 270  (478) & 13    & 16    & 34    & 6     & 311  (591) \\
      & CPC-SIS$_{(0.5,0.5)}$ & 5     & 5     & 6     & 1     & 110  (404) & 4     & 4     & 4     & 1     & 24   (371) \\
      & CPC-SIS$_{(0.75,0.75)}$ & 11    & 13    & 15    & 3     & 225  (453) & 13    & 39    & 14    & 5     & 299  (537) \\
\hline\hline
\end{tabular*}%
}
\end{table}%

\begin{table}
{ \footnotesize  
\tabcolsep 0pt \vspace*{-12pt}
\def\temptablewidth{1.0\textwidth}
\centering
\caption{\small Simulation results for Example 5, where $R_j$ indicates the median of the rank of the relevant predictors and MMS stands for the median of the minimum model size and its robust standard deviations (RSD) are given in parenthesis.  }  \label{tab5}%
\begin{tabular*}{\temptablewidth}{@{\extracolsep{\fill}} ll ccccc ccccc}
\hline\hline
      &       & \multicolumn{5}{c}{$\varepsilon\sim N(0,1)$}                    & \multicolumn{5}{c}{$\varepsilon\sim \frac{1}{3}Cauchy(0,1)$}   \\ \cline{3-7}\cline{8-12}
 $\rho$ &  Method$(\tau, \iota)$     & $R_1$    & $R_2$    & $R_3$    & $R_4$    & MMS (RSD) & $R_1$    & $R_2$    & $R_3$    & $R_4$    & MMS (RSD) \\  \hline
0.5   & SIS   & 12    & 12    & 377   & 4     & 488  (518) & 14    & 14    & 396   & 5     & 532  (494) \\
      & SIRS  & 192   & 215   & 910   & 188   & 910  (96)  & 236   & 222   & 938   & 193   & 938  (98) \\
      & DC-SIS & 336   & 305   & 511   & 320   & 744  (162) & 305   & 272   & 546   & 324   & 746  (154) \\
      & Kendall-SIS & 2     & 2     & 997   & 1     & 997  (15)  & 2     & 2     & 998   & 1     & 998  (10) \\
      & CC-SIS$_{(0.25,0.25)}$ & 3     & 3     & 841   & 2     & 841  (235) & 5     & 3     & 830   & 2     & 830  (263) \\
      & CC-SIS$_{(0.5,0.5)}$ & 3     & 4     & 886   & 2     & 886  (175) & 3     & 3     & 895   & 2     & 895  (207) \\
      & CC-SIS$_{(0.75,0.75)}$ & 3     & 3     & 817   & 2     & 817  (240) & 4     & 5     & 865   & 2     & 865  (184) \\
      & QC-SIS$_{(0.25)}$ & 183   & 249   & 713   & 208   & 796  (172) & 178   & 148   & 748   & 212   & 806  (170) \\
      & QC-SIS$_{(0.5)}$ & 269   & 241   & 672   & 307   & 823  (160) & 232   & 231   & 687   & 296   & 824  (172) \\
      & QC-SIS$_{(0.75)}$ & 223   & 174   & 701   & 259   & 829  (179) & 154   & 191   & 731   & 209   & 822  (174) \\
      & QPC-SIS$_{(0.25)}$ & 6     & 7     & 3     & 3     & 107  (167) & 7     & 9     & 3     & 4     & 109  (187) \\
      & QPC-SIS$_{(0.5)}$ & 4     & 5     & 2     & 3     & 53   (95)  & 5     & 5     & 3     & 5     & 77   (118) \\
      & QPC-SIS$_{(0.75)}$ & 5     & 7     & 3     & 3     & 75   (148) & 8     & 6     & 3     & 4     & 94   (165) \\
      & CPC-SIS$_{(0.25,0.25)}$ & 5     & 4     & 5     & 1     & 28   (58)  & 7     & 8     & 8     & 2     & 62   (110) \\
      & CPC-SIS$_{(0.5,0.5)}$ & 5     & 6     & 1     & 2     & 14   (27)  & 6     & 6     & 1     & 2     & 19   (34) \\
      & CPC-SIS$_{(0.75,0.75)}$ & 5     & 8     & 6     & 1     & 41   (88)  & 7     & 9     & 7     & 2     & 57   (75) \\   \hline
0.95  & SIS   & 5     & 9     & 538   & 4     & 581  (378) & 17    & 11    & 525   & 4     & 586  (445) \\
      & SIRS  & 190   & 231   & 894   & 199   & 894  (95)  & 246   & 220   & 895   & 186   & 895  (103) \\
      & DC-SIS & 441   & 303   & 508   & 239   & 771  (171) & 189   & 273   & 558   & 348   & 776  (174) \\
      & Kendall-SIS & 2     & 3     & 991   & 1     & 991  (32)  & 2     & 2     & 990   & 1     & 990  (46) \\
      & CC-SIS$_{(0.25,0.25)}$ & 3     & 4     & 806   & 2     & 806  (247) & 4     & 3     & 825   & 2     & 825  (254) \\
      & CC-SIS$_{(0.5,0.5)}$ & 3     & 4     & 831   & 2     & 831  (250) & 3     & 4     & 836   & 2     & 836  (257) \\
      & CC-SIS$_{(0.75,0.75)}$ & 3     & 4     & 804   & 2     & 804  (314) & 4     & 5     & 758   & 2     & 758  (315) \\
      & QC-SIS$_{(0.25)}$ & 326   & 225   & 624   & 157   & 795  (188) & 95    & 202   & 682   & 252   & 805  (192) \\
      & QC-SIS$_{(0.5)}$ & 400   & 244   & 597   & 226   & 804  (189) & 135   & 235   & 661   & 318   & 812  (193) \\
      & QC-SIS$_{(0.75)}$ & 312   & 168   & 650   & 143   & 785  (185) & 73    & 154   & 695   & 272   & 815  (178) \\
      & QPC-SIS$_{(0.25)}$ & 5     & 6     & 3     & 3     & 73   (165) & 12    & 6     & 3     & 4     & 128  (214) \\
      & QPC-SIS$_{(0.5)}$ & 4     & 5     & 3     & 3     & 47   (112) & 4     & 8     & 3     & 3     & 78   (115) \\
      & QPC-SIS$_{(0.75)}$ & 5     & 5     & 2     & 3     & 52   (146) & 8     & 6     & 3     & 4     & 105  (182) \\
      & CPC-SIS$_{(0.25,0.25)}$ & 6     & 5     & 6     & 1     & 31   (47)  & 6     & 5     & 7     & 2     & 55   (90) \\
      & CPC-SIS$_{(0.5,0.5)}$ & 5     & 6     & 1     & 2     & 20   (40)  & 6     & 5     & 1     & 2     & 18   (42) \\
      & CPC-SIS$_{(0.75,0.75)}$ & 6     & 6     & 5     & 1     & 41   (80)  & 6     & 7     & 7     & 2     & 46   (89) \\
\hline\hline
\end{tabular*}%
}
\end{table}%

\begin{table} 
{ \tiny 
\tabcolsep 0pt \vspace*{-12pt}
\def\temptablewidth{1.0\textwidth}
\centering
\caption{\small The overlaps of selected genes using various approaches for the rats data, where the screening threshold parameter is set as $\lfloor n/\log n \rfloor=25$ for each method and the CPC-SIS applies the algorithm in Case 1. }  \label{tab:d2}%
\begin{tabular*}{\temptablewidth}{@{\extracolsep{\fill}} ccccccccccccccccc}
\hline\hline
     &       &       &       &       & \multicolumn{3}{c}{QC-SIS$(\tau)$} & \multicolumn{3}{c}{QPC-SIS$(\tau)$} & \multicolumn{3}{c}{CC-SIS$(\tau,\iota)$} & \multicolumn{3}{c}{CPC-SIS$(\tau,\iota)$} \\
\cline{6-8}\cline{9-11}\cline{12-14}\cline{15-17}
 & SIS   & SIRS  & DC-SIS & Kendall-SIS & 0.25 & 0.5 & 0.75 & 0.25 & 0.5 & 0.75 & (0.25,0.25) & (0.5,0.5) & (0.75,0.75) & (0.25,0.25) & (0.5,0.5) & (0.75,0.75) \\
\hline
SIS   & 25    & 0     & 1     & 3     & 1     & 1     & 0     & 0     & 0     & 0     & 1     & 3     & 3     & 0     & 0     & 0 \\
SIRS  & 0     & 25    & 0     & 0     & 0     & 0     & 0     & 0     & 0     & 0     & 0     & 0     & 0     & 0     & 0     & 0 \\
DC-SIS & 1     & 0     & 25    & 1     & 2     & 2     & 1     & 0     & 0     & 0     & 0     & 0     & 1     & 0     & 0     & 0 \\
Kendall & 3     & 0     & 1     & 25    & 5     & 12    & 3     & 0     & 0     & 0     & 2     & 5     & 3     & 0     & 0     & 0 \\
QC-SIS$_{(0.25)}$ & 1     & 0     & 2     & 5     & 25    & 5     & 0     & 0     & 0     & 0     & 3     & 2     & 1     & 0     & 0     & 0 \\
QC-SIS$_{(0.5)}$ & 1     & 0     & 2     & 12    & 5     & 25    & 3     & 0     & 0     & 0     & 1     & 7     & 1     & 0     & 0     & 0 \\
QC-SIS$_{(0.75)}$ & 0     & 0     & 1     & 3     & 0     & 3     & 25    & 0     & 0     & 0     & 0     & 1     & 1     & 0     & 0     & 0 \\
QPC-SIS$_{(0.25)}$ & 0     & 0     & 0     & 0     & 0     & 0     & 0     & 25    & 3     & 2     & 0     & 0     & 1     & 2     & 1     & 0 \\
QPC-SIS$_{(0.5)}$ & 0     & 0     & 0     & 0     & 0     & 0     & 0     & 3     & 25    & 0     & 0     & 0     & 0     & 1     & 1     & 0 \\
QPC-SIS$_{(0.75)}$ & 0     & 0     & 0     & 0     & 0     & 0     & 0     & 2     & 0     & 25    & 0     & 0     & 0     & 0     & 1     & 0 \\
CC-SIS$_{(0.25,0.25)}$ & 1     & 0     & 0     & 2     & 3     & 1     & 0     & 0     & 0     & 0     & 25    & 1     & 0     & 0     & 0     & 0 \\
CC-SIS$_{(0.5,0.5)}$ & 3     & 0     & 0     & 5     & 2     & 7     & 1     & 0     & 0     & 0     & 1     & 25    & 1     & 0     & 0     & 0 \\
CC-SIS$_{(0.75,0.75)}$ & 3     & 0     & 1     & 3     & 1     & 1     & 1     & 1     & 0     & 0     & 0     & 1     & 25    & 0     & 0     & 0 \\
CPC-SIS$_{(0.25,0.25)}$ & 0     & 0     & 0     & 0     & 0     & 0     & 0     & 2     & 1     & 0     & 0     & 0     & 0     & 25    & 0     & 0 \\
CPC-SIS$_{(0.5,0.5)}$ & 0     & 0     & 0     & 0     & 0     & 0     & 0     & 1     & 1     & 1     & 0     & 0     & 0     & 0     & 25    & 1 \\
CPC-SIS$_{(0.75,0.75)}$ & 0     & 0     & 0     & 0     & 0     & 0     & 0     & 0     & 0     & 0     & 0     & 0     & 0     & 0     & 1     & 25 \\
\hline\hline
\end{tabular*}%
}
\end{table}

\begin{table}
{  \tiny
\tabcolsep 0pt \vspace*{-12pt}
\def\temptablewidth{1.0\textwidth}
\centering
\caption{\small Summary of top 10 gene probes selected by different screening methods for the rats data. ID means the selected gene ID and $p$-values are computed as $2(1-\Phi(|\sqrt{n} \widehat{\Omega}_1^{-1/2}\widehat{\varrho}_{Y,X}(0.5,0.5)|))$, where $\Phi$ is the cumulative distribution function of standard normal random variable. PE is the mean of prediction errors over 500 random partitions with the partition ratio of training sample to testing sample being $4:1$, where prediction error is defined as the average of $\{(Y_i-\widehat{Y}_i)^2, i\in \text{testing set}\}$.  PE1 and PE2 indicate that $\widehat{Y}_i$ is the predicted value via fitting a median regression model and linear model, respectively, using the top 10 genes selected. }  \label{tab:d3}%
\begin{tabular*}{\temptablewidth}{@{\extracolsep{\fill}} r cccc cccc cccc cccc} %
\hline\hline
      & \multicolumn{2}{c}{SIS} & \multicolumn{2}{c}{SIRS} & \multicolumn{2}{c}{DC-SIS} & \multicolumn{2}{c}{Kendall-SIS} & \multicolumn{2}{c}{QC-SIS(0.5)} & \multicolumn{2}{c}{CC-SIS(0.5,0.5)} & \multicolumn{2}{c}{QPC-SIS(0.5)} & \multicolumn{2}{c}{CPC-SIS(0.5, 0.5)} \\  \cline{2-3}\cline{4-5}\cline{6-7}\cline{8-9}\cline{10-11}\cline{12-13} \cline{14-15}\cline{16-17}
 Rank  & ID    & $p$-value & ID    & $p$-value & ID    & $p$-value & ID    & $p$-value & ID    & $p$-value & ID    & $p$-value & ID    & $p$-value & ID    & $p$-value \\
\hline
1     & 14770 & 2.1E-05 & 2828  & 1.000 & 146   & 6.4E-04 & 6083  & 8.3E-10 & 22641 & 5.2E-11 & 14726 & 4.4E-16 & 18602 & 0.473 & 1621  & 0.469 \\
2     & 21977 & 1.8E-05 & 20503 & 0.480 & 260   & 7.2E-04 & 5002  & 4.5E-10 & 14810 & 6.4E-12 & 6889  & 6.6E-14 & 4101  & 0.003 & 11288 & 0.001 \\
3     & 6436  & 2.0E-08 & 233   & 0.152 & 30768 & 2.3E-06 & 14726 & 4.4E-16 & 22339 & 8.0E-13 & 14701 & 6.2E-14 & 12365 & 0.141 & 12480 & 1.000 \\
4     & 4797  & 1.7E-08 & 3962  & 0.716 & 30745 & 1.6E-05 & 14810 & 6.4E-12 & 5002  & 4.5E-10 & 20898 & 7.2E-14 & 8399  & 1.000 & 4398  & 0.271 \\
5     & 21150 & 2.3E-07 & 7656  & 0.063 & 285   & 1.2E-04 & 25297 & 1.5E-11 & 20898 & 7.2E-14 & 22339 & 8.0E-13 & 5063  & 0.026 & 29604 & 0.467 \\
6     & 25573 & 4.5E-10 & 20453 & 0.468 & 30791 & 1.6E-07 & 5259  & 6.5E-12 & 31008 & 1.5E-07 & 23278 & 6.2E-13 & 9223  & 0.467 & 22679 & 1.000 \\
7     & 12127 & 9.4E-09 & 22023 & 0.047 & 3849  & 1.1E-04 & 5223  & 2.9E-10 & 26828 & 1.8E-08 & 25117 & 8.8E-14 & 21746 & 0.716 & 22267 & 0.065 \\
8     & 9235  & 1.6E-07 & 157   & 0.208 & 4626  & 1.1E-04 & 31008 & 1.5E-07 & 24529 & 4.2E-10 & 30548 & 6.9E-14 & 14019 & 0.717 & 17039 & 0.148 \\
9     & 3682  & 2.5E-06 & 2575  & 0.153 & 4490  & 2.8E-10 & 22339 & 8.0E-13 & 14414 & 1.5E-07 & 4512  & 8.0E-12 & 30361 & 0.277 & 11796 & 0.720 \\
10    & 8670  & 4.5E-11 & 2841  & 0.284 & 3967  & 2.1E-07 & 6021  & 1.7E-07 & 20724 & 2.8E-10 & 4712  & 8.4E-12 & 24759 & 0.010 & 20967 & 0.026 \\
\hline
PE1    & \multicolumn{2}{c}{0.0394} & \multicolumn{2}{c}{0.0252} & \multicolumn{2}{c}{0.0290} & \multicolumn{2}{c}{0.0310} & \multicolumn{2}{c}{0.0283} & \multicolumn{2}{c}{0.0330} & \multicolumn{2}{c}{0.0269} & \multicolumn{2}{c}{0.0247} \\
\hline
PE2   & \multicolumn{2}{c}{0.0377} & \multicolumn{2}{c}{0.0269} & \multicolumn{2}{c}{0.0349} & \multicolumn{2}{c}{0.0342} & \multicolumn{2}{c}{0.0344} & \multicolumn{2}{c}{0.0360} & \multicolumn{2}{c}{0.0307} & \multicolumn{2}{c}{0.0257} \\
\hline\hline
\end{tabular*}%
}
\end{table}%

\begin{table} 
{ \tiny 
\tabcolsep 0pt \vspace*{-12pt}
\def\temptablewidth{1.0\textwidth}
\centering
\caption{\small The overlaps of selected genes probes using various approaches for the breast cancer data, where the screening threshold parameter is set as $d_n= \lfloor n/\log n \rfloor=21$ for each method. The CPC-SISa$_1$ means the CPC-SIS in Case 1, the CPC-SISa$_2$ indicates the CPC-SIS in Case 2 and the CPC-SISa$_3$ stands for the CPC-SIS in Case 3.  }  \label{tab:d4}%
\begin{tabular*}{\temptablewidth}{@{\extracolsep{\fill}} ccccccccccccccccc}
\hline\hline
      &       &       &       &       & \multicolumn{3}{c}{QC-SIS$(\tau)$} & \multicolumn{3}{c}{QPC-SIS$(\tau)$} & \multicolumn{3}{c}{CC-SIS$(\tau,\iota)$} & CPC-SISa$_1$ & CPC-SISa$_2$ & CPC-SISa$_3$  \\   \cline{6-8}\cline{9-11}\cline{12-14}
    & SIS & SIRS & DC-SIS & Kendall-SIS & 0.25  & 0.5   & 0.75  & 0.25  & 0.5   & 0.75  & (0.25,0.25) & (0.5,0.5) & (0.75,0.75) & (0.5,0.5) & (0.5,0.5) & (0.5,0.5)  \\
\hline
SIS   & 21    & 8     & 7     & 4     & 1     & 7     & 5     & 1     & 1     & 0     & 2     & 1     & 1     & 0     & 1     & 0 \\
SIRS  & 8     & 21    & 13    & 9     & 1     & 15    & 2     & 1     & 0     & 0     & 3     & 5     & 0     & 0     & 1     & 0 \\
DC-SIS & 7     & 13    & 21    & 11    & 0     & 8     & 3     & 0     & 0     & 0     & 1     & 6     & 0     & 0     & 1     & 0 \\
Kendall & 4     & 9     & 11    & 21    & 0     & 3     & 4     & 0     & 0     & 0     & 0     & 3     & 0     & 0     & 0     & 0 \\
QC-SIS$_{(0.25)}$ & 1     & 1     & 0     & 0     & 21    & 1     & 0     & 0     & 0     & 0     & 1     & 0     & 0     & 0     & 0     & 0 \\
QC-SIS$_{(0.5)}$ & 7     & 15    & 8     & 3     & 1     & 21    & 0     & 1     & 0     & 0     & 3     & 7     & 0     & 0     & 1     & 0 \\
QC-SIS$_{(0.75)}$ & 5     & 2     & 3     & 4     & 0     & 0     & 21    & 0     & 1     & 0     & 0     & 0     & 5     & 0     & 0     & 0 \\
QPC-SIS$_{(0.25)}$ & 1     & 1     & 0     & 0     & 0     & 1     & 0     & 21    & 0     & 0     & 0     & 0     & 0     & 0     & 0     & 0 \\
QPC-SIS$_{(0.5)}$ & 1     & 0     & 0     & 0     & 0     & 0     & 1     & 0     & 21    & 0     & 0     & 0     & 0     & 0     & 0     & 0 \\
QPC-SIS$_{(0.75)}$ & 0     & 0     & 0     & 0     & 0     & 0     & 0     & 0     & 0     & 21    & 0     & 0     & 0     & 0     & 0     & 0 \\
CC-SIS$_{(0.25,0.25)}$ & 2     & 3     & 1     & 0     & 1     & 3     & 0     & 0     & 0     & 0     & 21    & 1     & 0     & 0     & 0     & 0 \\
CC-SIS$_{(0.5,0.5)}$ & 1     & 5     & 6     & 3     & 0     & 7     & 0     & 0     & 0     & 0     & 1     & 21    & 0     & 0     & 1     & 0 \\
CC-SIS$_{(0.75,0.75)}$ & 1     & 0     & 0     & 0     & 0     & 0     & 5     & 0     & 0     & 0     & 0     & 0     & 21    & 0     & 1     & 0 \\
CPC-SISa$_1$$_{(0.5,0.5)}$ & 0     & 0     & 0     & 0     & 0     & 0     & 0     & 0     & 0     & 0     & 0     & 0     & 0     & 21    & 0     & 0 \\
CPC-SISa$_2$$_{(0.5,0.5)}$ & 1     & 1     & 1     & 0     & 0     & 1     & 0     & 0     & 0     & 0     & 0     & 1     & 1     & 0     & 21    & 0 \\
CPC-SISa$_3$$_{(0.5,0.5)}$ & 0     & 0     & 0     & 0     & 0     & 0     & 0     & 0     & 0     & 0     & 0     & 0     & 0     & 0     & 0     & 21 \\
\hline\hline
\end{tabular*}%
}
\end{table} 

\begin{table}
{ \tiny
\tabcolsep 0pt \vspace*{-12pt}
\def\temptablewidth{1.0\textwidth}
\centering
\caption{\small Summary of top 10 gene probes selected by different screening methods for the breast cancer data. ID means the selected gene ID and $p$-values are computed as $2(1-\Phi(|\sqrt{n} \widehat{\Omega}_1^{-1/2}\widehat{\varrho}_{Y,X}(0.5,0.5)|))$, where $\Phi$ is the cumulative distribution function of standard normal random variable. PE is the mean of prediction errors over 500 random partitions with the partition ratio of training sample to testing sample being $4:1$, where prediction error is defined as the average of $\{(Y_i-\widehat{Y}_i)^2, i\in \text{testing set}\}$. PE1 and PE2 indicate that $\widehat{Y}_i$ is the predicted value via fitting a median regression model and linear model, respectively, using the top 10 genes selected.}  \label{tab:d5}
\begin{tabular*}{\temptablewidth}{@{\extracolsep{\fill}} r cccc cccc cc}
\hline\hline
      & \multicolumn{2}{c}{SIS} & \multicolumn{2}{c}{SIRS} & \multicolumn{2}{c}{DC-SIS} & \multicolumn{2}{c}{Kendall-SIS} & \multicolumn{2}{c}{QC-SIS(0.5)} \\
\cline{2-3}\cline{4-5}\cline{6-7}\cline{8-9}\cline{10-11}
Rank  & ID    & $p$-value & ID    & $p$-value & ID    & $p$-value & ID    & $p$-value & ID    & $p$-value \\
\hline
1     & 24032 & 3.0E-06 & 24032 & 0.000 & 8349  & 3.1E-05 & 17679 & 6.2E-02 & 24032 & 3.0E-06 \\
2     & 11913 & 1.2E-07 & 6841  & 0.000 & 24032 & 3.0E-06 & 20238 & 2.3E-02 & 22705 & 2.8E-02 \\
3     & 11870 & 2.9E-06 & 9164  & 0.001 & 13025 & 6.7E-04 & 10408 & 1.9E-03 & 6841  & 6.8E-08 \\
4     & 17439 & 6.9E-06 & 13025 & 0.001 & 23670 & 5.9E-03 & 1644  & 6.9E-06 & 14466 & 1.8E-03 \\
5     & 6841  & 6.8E-08 & 2172  & 0.013 & 20121 & 1.2E-07 & 8339  & 2.6E-03 & 4767  & 1.2E-05 \\
6     & 20938 & 2.3E-02 & 17439 & 0.000 & 6841  & 6.8E-08 & 14028 & 8.9E-05 & 5644  & 2.2E-04 \\
7     & 10692 & 2.0E-01 & 20121 & 0.000 & 15674 & 1.2E-06 & 23670 & 5.9E-03 & 20121 & 1.2E-07 \\
8     & 19897 & 1.5E-03 & 11870 & 0.000 & 1644  & 6.9E-06 & 12305 & 7.3E-04 & 23670 & 5.9E-03 \\
9     & 9164  & 1.5E-03 & 22705 & 0.028 & 5644  & 2.2E-04 & 3929  & 1.8E-03 & 13742 & 1.5E-05 \\
10    & 17050 & 2.2E-02 & 10408 & 0.002 & 20238 & 2.3E-02 & 14466 & 1.8E-03 & 17439 & 6.9E-06 \\
\hline
PE1    & \multicolumn{2}{c}{1.566} & \multicolumn{2}{c}{1.483} & \multicolumn{2}{c}{1.419} & \multicolumn{2}{c}{1.399} & \multicolumn{2}{c}{1.409} \\
\hline
PE2   & \multicolumn{2}{c}{1.550} & \multicolumn{2}{c}{1.398} & \multicolumn{2}{c}{1.378} & \multicolumn{2}{c}{1.366} & \multicolumn{2}{c}{1.367} \\
\hline
\hline
      & \multicolumn{2}{c}{CC-SIS(0.5,0.5)} & \multicolumn{2}{c}{QPC-SIS(0.5)} & \multicolumn{2}{c}{CPC-SISa$_1$(0.5, 0.5)} & \multicolumn{2}{c}{CPC-SISa$_2$(0.5, 0.5)} & \multicolumn{2}{c}{CPC-SISa$_3$(0.5, 0.5)} \\
\cline{2-3}\cline{4-5}\cline{6-7}\cline{8-9}\cline{10-11}
Rank  & ID    & $p$-value & ID    & $p$-value & ID    & $p$-value & ID    & $p$-value & ID    & $p$-value\\
\hline
1     & 12801 & 1.5E-06 & 11696 & 0.001 & 301   & 0.005 & 20121 & 0.000 & 4132  & 0.136 \\
2     & 13742 & 1.5E-05 & 672   & 0.005 & 18678 & 0.000 & 4356  & 0.001 & 17568 & 0.620 \\
3     & 402   & 6.2E-05 & 21944 & 0.021 & 3524  & 0.603 & 13084 & 0.008 & 5459  & 0.482 \\
4     & 4862  & 3.4E-04 & 6466  & 0.024 & 5422  & 0.021 & 13191 & 0.035 & 1079  & 0.352 \\
5     & 8349  & 3.1E-05 & 518   & 0.758 & 14782 & 0.023 & 6436  & 0.299 & 23942 & 0.002 \\
6     & 9158  & 1.9E-03 & 12635 & 0.022 & 21431 & 0.922 & 10179 & 0.192 & 14    & 0.922 \\
7     & 12074 & 6.8E-06 & 12567 & 0.609 & 5239  & 0.295 & 20102 & 0.185 & 1847  & 0.169 \\
8     & 14466 & 1.8E-03 & 7160  & 0.483 & 777   & 0.352 & 1299  & 0.179 & 3392  & 0.505 \\
9     & 18903 & 2.5E-02 & 21188 & 0.007 & 20958 & 0.132 & 1830  & 0.381 & 20369 & 0.367 \\
10    & 19774 & 8.1E-06 & 11916 & 0.495 & 4849  & 0.460 & 6025  & 0.467 & 390   & 0.920 \\
\hline
PE1    & \multicolumn{2}{c}{1.404} & \multicolumn{2}{c}{1.466} & \multicolumn{2}{c}{1.399} & \multicolumn{2}{c}{1.436} & \multicolumn{2}{c}{1.372} \\
\hline
PE2    &\multicolumn{2}{c}{1.290} & \multicolumn{2}{c}{1.437} & \multicolumn{2}{c}{1.345} & \multicolumn{2}{c}{1.304} & \multicolumn{2}{c}{1.289} \\
\hline\hline
\end{tabular*}%
}
\end{table}%

\end{document}